\input amstex
 \documentstyle{amsppt}

\nologo
 \pageno=1 \loadbold
  \leftheadtext{{\smc V. Mastropietro and M. Procesi}}
   \rightheadtext{{\smc  Lindstedt series for beam equations}}

\hsize=5 true in
 \vsize=8.2 true in
  \hoffset=.75 in
   \TagsOnRight
    \NoBlackBoxes

      \font\sevenrm=cmr7

\def\cal#1{\fam2 #1}

\topinsert
 \vskip-1.5
  \baselineskip {\vbox{\sevenrm\baselineskip 7pt

\noindent COMMUNICATIONS ON
 \hfill Website: http://AIMsciences.org\break
 \noindent PURE AND APPLIED ANALYSIS\hfil\break
  \line {Volume {\sevenbf ?}, Number {\sevenbf ?},
         Month {\sevenbf 2003}\hfil \eightpoint pp. ??--??}}
  }\endinsert

\footnote""{AMS 2000 Subject Classification:
   Primary 35B10, 35B32, 35L70, 47H15, 47N20, 58F27, 58F39.} 
\footnote""{Key words and phrases:
   Nonlinear wave equation, periodic solutions,
   Lindstedt series method, tree formalism, perturbation theory,
   Diophantine and irrationality conditions, Dirichlet boundary conditions.}
+
\bigskip

\document

\vglue 1\baselineskip

\centerline{\bf LINDSTEDT SERIES FOR 
PERIODIC SOLUTIONS }
\vskip.3cm
\centerline {\bf OF BEAM EQUATIONS WITH QUADRATIC AND }
\vskip.3cm
\centerline {\bf VELOCITY DEPENDENT
 NONLINEARITIES}

\bigskip

\medskip

\centerline{\smc Vieri Mastropietro }
\smallskip

{\eightpoint \centerline{ Dipartimento di Matematica}
\centerline{Universit\`a di Roma ``Tor Vergata'', Roma, I-00133}
 }
\medskip
\centerline{\smc   Michela Procesi}

\smallskip

{\eightpoint
 \centerline{ Dipartimento di Matematica}
  \centerline{  Universit\`a di Roma Tre, Roma, I-00146}}

\medskip

 \centerline{(Communicated by Vieri Mastropietro)}

 \bigskip
{\eightpoint{\narrower\smallskip\noindent
 {\bf Abstract.}
{ We prove the existence of small amplitude periodic solutions,
for a large Lebesgue measure set of frequencies,
in the nonlinear beam equation with a weak quadratic 
and velocity dependent nonlinearity and with
Dirichlet boundary conditions. Such nonlinear PDE
can be regarded as a simple model describing oscillations of flexible
structures like suspension bridges in presence of an uniform wind flow.  
The periodic solutions are explicitly constructed by means of 
a perturbative expansion 
which can be considered
the analogue of the Lindstedt series expansion
for the invariant tori in classical mechanics. 
The periodic solutions are not analytic but defined only in a
Cantor set, and resummation techniques of divergent powers series
are used in order to control the small divisors problem.
}
\smallskip}}

\medskip



\font\tenmib=cmmib10
\font\sevenmib=cmmib10 scaled 800

 2

\font\sc=cmcsc10

\font\ottorm=cmr8
\textfont5=\tenmib\scriptfont5=\sevenmib\scriptscriptfont5=\fivei

\font\msytw=msbm10 scaled\magstep1

\font\msytwww=msbm7 scaled\magstep1
\font\indbf=cmbx10 scaled\magstep2

\font\ottorm=cmr8\font\ottoi=cmmi8\font\ottosy=cmsy8%
\font\ottobf=cmbx8\font\ottott=cmtt8%
\font\ottocss=cmcsc8%
\font\ottosl=cmsl8\font\ottoit=cmti8%
\font\sixrm=cmr6\font\sixbf=cmbx6\font\sixi=cmmi6\font\sixsy=cmsy6%
\font\fiverm=cmr5\font\fivesy=cmsy5\font\fivei=cmmi5\font\fivebf=cmbx5%

\def\ottopunti{\def\rm{\fam0\ottorm}%
\textfont0=\ottorm\scriptfont0=\sixrm\scriptscriptfont0=\fiverm%
\textfont1=\ottoi\scriptfont1=\sixi\scriptscriptfont1=\fivei%
\textfont2=\ottosy\scriptfont2=\sixsy\scriptscriptfont2=\fivesy%
\textfont3=\tenex\scriptfont3=\tenex\scriptscriptfont3=\tenex%
\textfont4=\ottocss\scriptfont4=\sc\scriptscriptfont4=\sc%
\textfont\itfam=\ottoit\def\it{\fam\itfam\ottoit}%
\textfont\slfam=\ottosl\def\sl{\fam\slfam\ottosl}%
\textfont\ttfam=\ottott\def\tt{\fam\ttfam\ottott}%
\textfont\bffam=\ottobf\scriptfont\bffam=\sixbf%
\scriptscriptfont\bffam=\fivebf\def\bf{\fam\bffam\ottobf}%
\setbox\strutbox=\hbox{\vrule height7pt depth2pt width0pt}%
\normalbaselineskip=9pt\let\sc=\sixrm\normalbaselines\rm}
%


\newcount\driver 
\newdimen\xshift \newdimen\xwidth
\def\ins#1#2#3{\vbox to0pt{\kern-#2 \hbox{\kern#1
#3}\vss}\nointerlineskip}

\def\insertplot#1#2#3#4#5{\par%
\xwidth=#1 \xshift=\hsize \advance\xshift by-\xwidth \divide\xshift by 2%
\yshift=#2 \divide\yshift by 2%
\line{\hskip\xshift \vbox to #2{\vfil%
\ifnum\driver=0 #3
\special{ps: plotfile #4.eps} 
\fi
\ifnum\driver=1 #3 \includegraphics{#4.eps}\fi
\ifnum\driver=2 #3
\ifnum\mgnf=0\special{#4.eps 1. 1. scale} \fi
\ifnum\mgnf=1\special{#4.eps 1.2 1.2 scale}\fi
\fi }\hfill \raise\yshift\hbox{#5}}}

\def\insertplotttt#1#2#3{\par%
\xwidth=#1 \xshift=\hsize \advance\xshift by-\xwidth \divide\xshift by 2%
\yshift=#2 \divide\yshift by 2%
\line{\hskip\xshift \vbox to #2{\vfil%
\includegraphics{#3.eps}}\hfill}}

\def\insertplott#1#2#3{\par%
\xwidth=#1 \xshift=\hsize \advance\xshift by-\xwidth \divide\xshift by 2%
\yshift=#2 \divide\yshift by 2%
\line{\hskip\xshift \vbox to #2{\vfil%
\ifnum\driver=0
\special{ps: plotfile #3.eps} 
 \fi
\ifnum\driver=1 \includegraphics{#3.eps}\fi
\ifnum\driver=2
\ifnum\mgnf=0\special{#3.eps 1. 1. scale} \fi 
\ifnum\mgnf=1\special{#3.eps 1.2 1.2 scale}\fi
\fi }\hfill}}

\newdimen\xshift \newdimen\xwidth \newdimen\yshift
\def\eqfig#1#2#3#4#5{
\par\xwidth=#1 \xshift=\hsize \advance\xshift
by-\xwidth \divide\xshift by 2
\yshift=#2 \divide\yshift by 2
\line{\hglue\xshift \vbox to #2{\vfil
\ifnum\driver=0 #3
\special{ps: plotfile #4.ps} 
\fi
\ifnum\driver=1 #3 \includegraphics{#4.ps}\fi
\ifnum\driver=2 #3 \special{
\ifnum\mgnf=0 #4.ps 1. 1. scale \fi
\ifnum\mgnf=1 #4.ps 1.2 1.2 scale\fi}
\fi}\hfill\raise\yshift\hbox{#5}}}


\let\a=\alpha   \let\g=\gamma  \let\d=\delta \let\e=\varepsilon
  \let\h=\eta   \let\th=\theta  
\let\m=\mu    \let\n=\nu         \let\p=\pi    
\let\s=\sigma \let\t=\tau    
   \let\o=\omega
   \let\Th=\Theta\let\L=\Lambda 
         
\let\O=\Omega 
\def\ka{{(k)}}
\def\\{\hfill\break} \let\==\equiv

\let\io=\infty 

\let\0=\noindent

\def\ie{\hbox{\it i.e.\ }}
\let\dpr=\partial

\def\tende#1{\,\vtop{\ialign{##\crcr\rightarrowfill\crcr
 \noalign{\kern-1pt\nointerlineskip}
 \hskip3.pt${\scriptstyle #1}$\hskip3.pt\crcr}}\,}
\def\otto{\,{\kern-1.truept\leftarrow\kern-5.truept\to\kern-1.truept}\,}

\def\R{{\cal R}}\def\D{{\cal D}}
 \def\VV{{\cal V}}
\def\CC{{\cal C}}
\def\NN{{\cal N}}\def\II{{\cal I}}
\def\RR{{\cal R}}\def\LL{{\cal L}} 
\def\GG{{\cal G}}

\def\T#1{{#1_{\kern-3pt\lower7pt\hbox{$\widetilde{}$}}\kern3pt}}
\def\VVV#1{{\underline #1}_{\kern-3pt
\lower7pt\hbox{$\widetilde{}$}}\kern3pt\,}
\def\W#1{#1_{\kern-3pt\lower7.5pt\hbox{$\widetilde{}$}}\kern2pt\,}

  \def\sign{{\text sign}\,}
\def\indica{\leaders \hbox to 0.5cm{\hss.\hss}\hfill}
\def\guida{\leaders\hbox to 1em{\hss.\hss}\hfill}


\mathchardef\aa   = "050B
\mathchardef\bb   = "050C
\mathchardef\xxx  = "0518
\mathchardef\hhh  = "0511
\mathchardef\zzzzz= "0510
\mathchardef\oo   = "0521
\mathchardef\lll  = "0515
\mathchardef\mm   = "0516
\mathchardef\Dp   = "0540
\mathchardef\H    = "0548
\mathchardef\FFF  = "0546
\mathchardef\ppp  = "0570
\mathchardef\nn   = "0517
\mathchardef\pps  = "0520
\mathchardef\XXX  = "0504
\mathchardef\FFF  = "0508
\mathchardef\nnnnn= "056E

\def\to{\rightarrow}

\def\eps{\varepsilon}

\let\ciao=\bye
\def\qed{\raise1pt\hbox{\vrule height5pt width5pt depth0pt}}

 \def\Val{{\text Val}}
\def\indic{\hbox{\raise-2pt \hbox{\indbf 1}}}

\def\RRR{\hbox{\msytw R}}

\def\NNN{\hbox{\msytw N}} 
 \def\ZZZ{\hbox{\msytw Z}}
 \def\zzz{\hbox{\msytwww Z}}
\def\TTT{\hbox{\msytw T}} 

\def\QQQ{\hbox{\msytw Q}}

\bigskip
\noindent {\bf 1. Introduction and  main Results}
\smallskip
\noindent {\bf 1.1)}
The search of periodic solutions in nonlinear wave equations
has attracted a wide interest in recent times.
In the finite dimensional case the problem has its analogous
in the study of periodic orbits close to elliptic equilibrium points:
results of existence have been obtained in such a case starting from  Lyapunov
[20].
Systems with infinitely many degrees of freedom (as
the nonlinear wave equation, the nonlinear Schr\"odinger equation and
other PDE systems) have been studied much more recently;
the problem is much more difficult because of the presence of a
{\it small divisors problem}, which is absent in the finite dimensional case,
and one has to prove an infinite dimensional KAM theorem to overcome such 
difficulty. Periodic or quasi periodic solutions in PDE
have been obtained for instance in [21],[9],[17],[3],[7],[6]
by a Lyapunov-Schmidt decomposition together with
KAM methods. Generally the nonlinear terms are assumed {\it odd}
and velocity-independent, as such features considerably
simplify the analysis.
A velocity dependent non linearity has  been considered in  [5], in which the
string equation with a nonlinear term $u_t^2$ and 
{\it periodic} boundary conditions is considered. 
The recent papers [3] and [4] consider the  massless string equation,
 under Dirichlet boundary conditions, with  velocity independent (but otherwise quite general) nonlinearities.

Aim of this paper is to construct periodic solutions in a {\it beam
equation} with an even and velocity dependent nonlinearity
and {\it Dirichelet} boundary conditions;

$$\left\{\eqalign{ &v_{tt} + \partial^4_x v +\m v=  a v^2+ b v_t^2  ,  \cr & v(0,t) = v(\p,t) = 0 ,  \cr}\right.\tag 1.1 $$
where $ a,b,\m$ are suitable parameters. As it will appear
clear in the following, our results could be easily extended
to include more general nonlinearities. With respect to [5], we have considered the beam
instead of the wave equation, leading to a simpler small divisor
problem; on the other hand 
Dirichelet boundary conditions and even 
nonlinearities introduces various regularity problems which are not 
present in the case of periodic boundary conditions considered in [5].

The interest of 1.1 lies moreover in the fact that it
can be regarded as a simple model describing oscillations of flexible
structures; for instance, see [16],[8], a suspension bridge   
subjected to elastic forces 
due to suspensions and to forces 
caused by a uniform wind-flow has been described by a beam equation 
with a nonlinear terms quadratic in 
$v$ (describing the anharmonic elastic forces)
and depending also from $v_t$ 
(to take into account the forces due to the wind flow).
Another applications of PDE 
with this kind of nonlinear terms is in 
[17] to describe the oscillations of the atmosphere on the flat earth.
In the literature there is no  proof of existence of periodic solutions
in a large set for such a problem. We will construct such solutions
generalizing to the present case the approach based on Lindstedt
series expansion already adopted first 
in [12] to prove the existence of periodic solutions in a zero measure set,
and later on generalized to construct periodic
solutions in a large measure set in [13],[14].

We call $\o_m=\sqrt{m^4+\m}$ so that if 
the non linear terms are absent $a=b=0$
every solution of 1.1 can be written as
$$ v(x,t) = 
\sum_{m=0}^\io A_{m}\cos(\o_m t+\theta_{m}) \sin {m x} ,
\tag 1.2 $$
where $\theta_{m}$ is an arbitrary phase.
In particular if $\m\notin \QQQ$ the only $\o_1$ periodc solutions are of 
the form:
$$ \pm \sqrt{\e} \cos \o_{1} t \sin x \tag 1.3 $$
for all values of the real parameter 
$\e$. We will show that  
also if the nonlinear term
is added to 1.1, that is $a\not=0$ or $b\not=0$,
 periodic solutions close to 1.3 exists.
\smallskip 
\noindent {\bf 1.2)} 
To face the small divisor problem, some Diophantine conditions must be imposed
on the mass $\m$.
\smallskip
\noindent{\bf Definition 1.} {\it 
We call $M(\g)$, $\g\le 2^{-6}$, the set $\m\in [0,\m_0]$, $\m_0={1\over 8}$ 
verifying the following Diophantine condition
$$|\o_1 n \pm \o_m|
\ge \g |n|^{-\t_0}
\qquad \forall n\in\ZZZ \setminus \{ 0 \} \hbox{ and }
\forall m \in \NNN\setminus\{1\}\tag 1.4 $$
$$| \o_1 n \pm \o_m\pm \o_{m'}|
\ge \g |n|^{-\t_0}
\qquad \forall n\in\ZZZ \setminus \{ 0 \} \hbox{ and }
\forall m,m' \in \NNN\setminus\{1\}$$}
\smallskip
It will be shown in Appendix A1
that the set of $\m$ verifying 1.4, for some positive $\g$, 
is of measure $O(\m_0)$ provided that $\t_0\ge 4$ and $\g$ is small enough.
\smallskip
Our main result is the following Theorem.
\smallskip
\0 {\bf Theorem 1.} {\it Generically in $a,b$, 
for any $\m\in M(\g)$ there exists an $\e_0>0$ and
a Cantor set $\CC(\g)\subset (0,\e_0)$ verifying
$\lim_{\e\to 0^+}{1\over \e} \hbox{ meas } (\CC(\g)
\cap (0,\e))=1$ such that
for all $\e\in \CC(\g)$ there exists  a periodic solution $v(x,\O t):\TTT^2\to \RRR$ of 1.1
, with $\O=\o_1+\e$, of the form
$$v(x,\O t,\e)=\sqrt\e u(x,\O t;\e) =\sqrt{\e}\sum_{n\in\zzz}\sum_{m=1}^\io
e^{in\O t} \sin (m x) u_{n,m}\;\tag 1.5 $$
with $u_{n,m}=u_{-n,m}$ and
$$|u_{n,m}|\le {C_0 e^{-\s |n|} \over m^7}\tag 1.6$$
with suitable constants $C_0,\s$.}
\smallskip
Note that in presence of odd nonlinearities, like $v^3$, one can continue
the periodic solution in an analytic solution both
is space and time, see [13]
; on the contrary, in presence of even  
or velocity depending nonlinearities, like in the present case, 
the periodic solutions are {\it not analytic} in space
and this lack of regularity is reflected in some complications
in their constructions.
\smallskip 
\noindent {\bf 1.3)} By inserting (1.6) in (1.1)
we get a closed equation for the 
coefficients $u_{n,m}(\e)\equiv u_{n,m}$
$$ u_{n,m} \left[ - \O^{2} n^2
+\o_m^2 \right] =\sqrt{\e} \hat f_{n,m}(u) .
\tag 1.7 $$
where
$$\hat f_{n,m}(u)={\O \over 2\pi^2}\int_0^\pi dx\int_0^{2\pi/\O}dt \sin(m x)
e^{-i n \O t}( a u^2+ b u_t^2)\tag 1.8$$

More explicitly, see Appendix A2, 1.8  can be written as
$$\hat f_{n,m}=\sum^*_{n_1,m_1\atop n_2,m_2}v_{m,m_1,m_2} \d_{n_1+n_2,n}(a 
+b (i\O n_1)
(i\O n_2))u_{n_1,m_1}u_{n_2,m_2} \tag 1.9$$
where $\sum^*$ means that the sum is over $m,m_1,m_2$ 
such that $m\pm m_1\pm m_2=odd$ and
$$v_{m,m_1,m_2}={ 4 m m_1 m_2\over \pi (m^2-(m_1-m_2)^2)(m^2-(m_1+m_2)^2)}\tag 1.10$$
One could try to write a power series
expansion in $\e$ for $u(x,t)$, using 1.7
to get recursive equations for the coefficients.
However by proceeding in this way one finds that
the coefficient of order $k$ is given by a sum of terms
some of which of order $O(k!^\a)$, for some constant $\a$.
This is the same phenomenon occurring
in the Lindstedt series for invariant KAM tori [10],[11]
in the case of quasi-integrable Hamiltonian systems;
in such a case however one can show that
there are {\it cancellations} between the terms
contributing to the coefficient of order $k$,
which at the end admits a bound $C^k$, for a suitable constant $C$.
On the contrary such cancellations are absent in the present case
and we have to proceed in a different way, 
essentially equivalent to a {\it resummation}.

We write
$$\left\{\eqalign{& \h u_{1,1} \equiv \h q= 
\hat f_{1,1}(u)\quad\quad\quad\quad\quad \hbox{ if }(|n|,m)=(1,1)\quad\quad\quad\qquad\qquad\qquad\quad\; (1.11)\;\cr & u_{n,m} \left[ - \O^{2} n^2
+\o_m^2 +n\n_{n,m}\right] \equiv g^{-1}_{n,m}u_{n,m} =\h( \hat f_{n,m}(u) +nl_{n,m}u_{n,m})\; \hbox{ otherwise }.\cr}\right.
$$
Naturally equations 1.11 coincides with 1.12 provided that:
$$\h=\sqrt{\e}\quad \n_{n,m}= \h l_{n,m}\tag 1.12$$
We introduce the following definition.
\smallskip
\noindent{\bf Definition 2.}
{\it We define $\cal D$, subset of $(\e,\n)\in \RRR^+\times l_\io$,
as 
$${\cal D}:= \{(\e,\n): 0<\e< \e_0\,,\;
\max_{n,m}|\n_{n,m}| < c \e_0,\;\n_{n,m}=0\;{\text if }
\,|\o_1 |n|-m^2|\ge 1+\e_0 |n| 
\}\tag 1.13$$ 
We define $\L$ the set of $(n,m)$ such that
$|\o_1 |n|-m^2|\le 1+\e_0 |n| $.

For any $\m\in M( 4\g)$, $\t>\t_0+5$, we define a subset ${\cal D}(\g)\subset \D$ of couples
$(\e,\n)\in \cal D$ verifying  the following Diophantine conditions

$$\left| \O n \pm \sqrt{\o_m^2+n\n_{n,m}} \right|
\ge \g |n|^{-\t}
\qquad \forall n\in\ZZZ \setminus \{ 0 \} \hbox{ and }
\forall m \in \NNN\setminus\{1\}\tag 1.14$$
$$\left| \O (n_2-n_1) \pm \sqrt{\o_{m_1}^2+n_1\n_{n_1,m_1}}\pm
\sqrt{\o_{m_2}^2+n_2\n_{n_2,m_2}}
\right| \ge \g |n_2-n_1|^{-\t}$$
$$\forall  n_1,n_2\in\ZZZ \setminus \{ 0 \} \hbox{ and }
\forall m_1\neq m_2 \in \NNN:\quad 
|\o_1 |n_i|-m_i^2|\le 1+\e_0 |n_i|,
i=1,2
\tag 1.15$$
%

\smallskip
}
We call 1.14 and 1.15, respectively the first and second Melnikov conditions.
Our strategy in order to prove Theorem 1
is the following: 

1) First 
we consider $(\e,\n)$ as independent parameters 
belonging to $D(\g)$ (so that the Melnikov conditions are verified)
and we show that it is possible to 
{\it find} a proper $l_{n,m}(\h,\e,\n)$, well defined for $|\h|\le \h_0$ 
and $(\e,\n)\in D(\g)$, such that 1.11 admits
a solution $u_{n,m}$ analytic in $\h$;
both $u_{n,m}$ and $l_{n,m}$
are expressed by convergent power series in 
$\h$. Using a technique inspired by [4], we extend $l_{n,m}$ to  a $C^1$ 
function,$l_{n,m}^E$,defined  on the 
square $\D$; $l_{n,m}^E$ coincieds with 
$l_{n,m}$ in the set $D(2\g)$.

2)The solution  $u_{n,m}$ defined above
is a solution of 1.7 
only if 1.12 is verified; we show
(Proposition 2) that we can find $\n=\n(\e)$ so that 
1.12 is verified for all $(\e,\n(\e))\in D(2\g)$; more precisely $\n(\e)$ 
solves the equation
 $\n_{n,m}=\sqrt{\e} l^E_{n,m}(\sqrt{\e},\e,\n)$: hence
replacing $\n_{n,m}$ with 
$\n_{n,m}(\e)$ in the expansion for $u_{n,m}$ we get the 
solution of 1.7.
\smallskip
\noindent{\bf Proposition 1}:
{\it Assume that $\m\in M(4\g)$ and $(\n,\e)\in D(\g)$. Let  $C_0,C_1,C_2,\s$
be positive constants. 
It is possible to find 
a sequence 
$$\{l_{n,m}(\h,\e,\n)\}_{(n,m)\in \ZZZ^2\setminus \{(\pm 1,1)\}}\tag 1.16$$  
such that:

(i) 
There exists a unique solution $u(\h,\n,\e;x,t)$, analytic in $t$ 
and $C^5$ in $x$, of equation 1.11; 
$u$ is analytic in $\h$ for $|\h|\leq \h_0$ and is such that:
$$|u(\h,\n,\e;x,t)-u_{1,1}(\n,\e)\cos \O t\sin x|\leq |\h| C_0.\tag 1.17 $$
for a proper $u_{1,1}(\n,\e,) $.

(ii) The sequence $l_{n,m}(\h,\e,\n)$ is analytic in $\h$ and  uniformly bounded for $(\e,\n)\in \D(\gamma)$:
$$ |l(\h,\e,\n)|_\io\equiv \max_{n,m}|l_{n,m}|\leq C_1|\h|.\tag 1.18$$ 

(iii) The functions $u_{n,m}(\h,\e,\n)$ and $l_{n,m}(\h,\e,\n)$  can be
 extended to  $C^1$ functions,
denoted by  $u_{n,m}^E(\h,\e,\n)$, $l_{n,m}^E(\h,\e,\n)$,
on the set $\D$, such that
$$ l_{n,m}^E(\h,\e,\n)=  l_{n,m}(\h,\e,\n)\quad \forall \; (\e,\n)\in \D(2\g) \tag 1.19$$
 The same is true for $u^E_{n,m}$.

(iv) $l_{n,m}^E(\h,\e,\n) $ respects the bounds:
$$|l^E(\h,\e,\n)|_\infty \leq |\h| C_2\,,
\quad |\partial_\e l^E(\h,\e,\n)|_\io 
\leq |\h| C_2 \,,\quad |\partial_{\n_{n,m}} l^E(\h,\e,\n)|_\infty \leq|\h| C_2 
\,, \quad\tag 1.20$$
$$|\sum_{(n,m)\in \L}\partial_{\n_{n,m}} l^E(\h,\e,\n)|_\infty \leq|\h| C_2\,,  \quad |u^E_{n,m}(\h,\e,\n)| \leq |\h| {1\over m^7} C_2 e^{-\s |n|}\tag 1.21$$}
\smallskip

Once we have proved Proposition 1, we solve the compatibility equation for the extended counterterm function $l_{n,m}^E(\h=\sqrt{\e},\e,\n) $ which is well defined provided that we choose $\e_0$ so that $\e_0<\h_0^2$.
\smallskip
\noindent{\bf Proposition 2.}{\it For all $(n,m)\neq (\pm 1,1)$, exist  $C^2$ functions 
$\n_{n,m}(\e):(0,\e_0)\to (-c\e_0,c\e_0)$ such that

(i) $\n_{n,m}(\e)$ verifies
$$\n_{n,m}(\e)= \sqrt{\e} l^E_{n,m}(\sqrt{\e},\e,\n_{n,m}(\e));\tag 1.22
$$ 
and is such that
 $$|\n_{n,m}(\e)|\leq C\e\,, \quad |\partial_\e \n_{n,m}(\e)|
\leq C\tag 1.23$$ 
for a suitable constant $C$;
\vskip.5cm
(ii) the set $\CC\equiv \CC(2\g)$ defined by $\e\in (0,\e_0)$ and the conditions:
$$\left|\O n - m^2\right|> 4\g |n|^{-\t_0}\tag 1.24 $$
$$\left| \O n \pm \sqrt{\o_m^2+n\n_{n,m}(\e)} \right|
\ge 2\g |n|^{-\t}
\qquad \forall n\in\ZZZ \setminus \{ 0 \} \hbox{ and }
\forall m \in \NNN\setminus\{1\}\tag 1.25$$
$$\left| \O (n_2-n_1) \pm (\sqrt{\o_{m_1}^2+n_1\n_{n_1,m_1}(\e)}
\pm \sqrt{\o_{m_2}^2+n_2\n_{n_2,m_2}(\e)})
 \right| \ge 2\g |n_2-n_1|^{-\t}\tag 1.26$$
$$\forall  
n,n_2\in\ZZZ \setminus \{ 0 \} \hbox{ and }
\forall m_1,m_2 \in \NNN\quad m_1-m_2\neq 0,
|\o_1 |n_i|-m_i^2|\le 1+\e_0 |n_i|,
i=1,2$$ 
has large relative Lebesgue measure, namely
$\lim_{\e\to 0^+}{1\over \e} \hbox{ meas } (\CC(\g)
\cap (0,\e))=1$.
}
\smallskip
\smallskip 
\0{\bf 1.4)}
Theorem 1 is an easy consequence of Proposition 1 and 2.
\smallskip
\0 {\it Proof of the Theorem 1.} 
We start by choosing  $\g$ and $\m\in M(4\g)$ and keep $\e_0$ as a parameter; by Proposition 1 (i) 
for all $(\e,\n)\in \D(\g)$ we can find a sequence $l_{n,m}$ so that there 
exists a unique solution $u(\h,\n,\e;x,t)$ of  1.11 
for all $|\h|\le \h_0$ where $\h_0$ depends only on $\g$ for $\e_0$ small enough.
 By Proposition 1 (iii) the sequence $l_{n,m}$ and the solution $u(\h,\n,\e;x,t)$ can be extended to $C^1$ functions ( denoted by $l^E,u^E$) for all $(\e,\n)\in \D$.
Moreover  $l^E_{n,m}(\e,\n)= l_{n,m}(\e,\n),$ $u^E_{n,m}(\e,\n)= u_{n,m}(\e,\n) $ for all $(\e,\n)\in \D(2\g)$.

Equation 1.11 coincides with our original equation 1.7 
provided that the compatibility equations 1.12 are satisfied. Now we fix $\e_0<\h_0^2$ so that  $l^E_{n,m}(\h=\sqrt{\e},\e,\n)$ and  $u^E_{n,m}(\h=\sqrt{\e},\e,\n)$ are well defined.
By Proposition 2 (i) there exists a sequence $\n_{n,m}(\e)$ which satisfies the extended compatibility equation 1.12.
Finally by Proposition 2(ii) the Cantor set $\CC(2\g)$ is well defined and of large relative measure.
 
For all $\eps\in \CC(2\g)$ 
we have that the couple $(\eps,\n(\e))$ is by definition in $\D(2\g)$ so  that
by Proposition 1(iii): $$\align &l_{n,m}(\sqrt{\e},\eps,\n(\e))=l^E_{n,m}(\sqrt{\e},\eps,\n(\e))\tag 1.27 
\\ &u(\sqrt{\e},\eps,\n(\e);x,t)=u^E(\sqrt{\e},\eps,\n(\e);x,t).\endalign $$ 
so that  $u(\sqrt{\e},\e,\n(\e);x,t)$ solves equation  1.11 for $\h=\sqrt{\e}$.
So by Proposition 2(i) $\n(\e)$ solves the true compatibility equation  1.12
$$\n_{n,m}(\e)=\sqrt{\e}l_{n,m}(\sqrt{\e},\e,\n(\e))\tag 1.28 $$ 
for all $\e\in \CC(2\g)$.
Then   $\sqrt{\e}u(\sqrt{\e},\e,\n(\e);x,t)$ 
is a true non trivial solution of our equation 1.1 in $\CC(2\g)$.
\qed
\smallskip
In the rest of the paper we prove Proposition 1 and 2. 

\medskip

\bigskip
\noindent {\bf 2. Lindstedt series and tree expansion.}

\smallskip 
\noindent{\bf 2.1)} In this section we find a formal solution $u_{n,m}$ of  1.11
as power series on $\h$; the solution $u_{n,m}$ is parameterized by the coefficients 
$l_{n,m}$ and it will be written in the form of a tree expansion.

We assume for  $l_{n,m}(\h,\e,\n), u_{n,m}(\h,\e,\n)$ with $(n,m)\neq (\pm 1,1)$, a formal series expansion in $\h$:
$$l_{n,m}(\h,\e,\n)= \sum_{k=2}^{\infty}\h^{k-1} l_{n,m}^\ka \,,\qquad u_{n,m}(\h,\e,\n)= \sum_{k=1}^{\infty}\h^k u_{n,m}^\ka\tag 2.1$$  
for all $(n,m)\neq(\pm 1,1)$.
By definition we set $q=u_{\pm 1,1}^{(0)}
$ and $u_{\pm 1,1}^\ka=0$. 
Inserting the series expansion in the second equation  1.11 we obtain the recursive equations:
$$  u_{n,m}^{(k)} = g_{n,m}\Big(n\sum_{r=2}^{k-1}l_{n,m}^{(r)} u^{(k-r)}_{n,m} +$$
$$\sum^*_{n_1,m_1\atop n_2,m_2}\sum_{k_1+k_2=k-1} 
\d_{n_1+n_2,n}(a-b \O^2 n_1n_2) v_{m,m_1,m_2}u_{n_1,m_1}^{(k_1)}u_{n_2,m_2}^{(k_2)}\Big)\equiv $$ $$g_{n,m} \Big(n\sum_{r=2}^{k-1}l_{n,m}^{(r)} u^{(k-r)}_{n,m} +F_{n,m}^{(k)}\Big)   
\tag 2.2$$
where $ g_{n,m}$ (defined in equation 1.11) is called the propagator. 
It holds the following Lemma.
\smallskip
\noindent{\bf Lemma 1.}{\it  For all $(n,m)\neq (\pm 1,1)$ we have that  $u_ {n,m}^\ka=0$ when $|n|>k+1$ or $m$ is even.}
\smallskip
{\it Proof.} We proceed by induction. By definition $F^{(0)}_{n,m}=0$ so that $u_{n,m}^{(0)}=0$ if $(n,m)\neq (\pm 1,1)$.

Now suppose that our claim holds for all $ (n,m)\neq (\pm 1,1)$ and $r<k$. Equations  2.2 are recursive 
so that $F^\ka_{n,m}  $ is a quadratic polynomial sum of monomials of the form   
$v(m_1,m_2,m)u_{n_1,m_1}^{(h_1)} u_{n_2,m_2}^{(h_2)} $ such that  the $m_i$ are odd, $|n_i|<h_i $, $n=n_1+n_2$ and $h_1+h_2=k-1$. This implies that $F^\ka_{n,m}  $ can be nonzero  only if  $n=n_1+n_2\leq h_1+h_2+2=k+1$. In the same way the linear terms $l^{(r)}_{n,m}u_{n,m}^{k-r-1} $ can be non zero only if $|n|\le k-1 $.

Finally the factor $v(m,m_1,m_2)=0$ if $m_1+m_2+m$ is even and by the inductive hypothesis $m_1$ and $m_2$ are odd so that $m$ must be odd as well. 
\qed
\smallskip 
We introduce a smooth partition of the unity in the following way.
Let $\chi(x)$ be a $C^\io$ non-increasing function such that
$\chi(x)=0$ if $|x|\le \g$ and $\chi(x)=1$ if $|x|\ge 2\g$;
moreover $|\chi'(x)|\le \g^{-1}$.
Let $\chi_h(x)=\chi(2^{h}x)-\chi(2^{h+1}x)$ for
$h \ge 0$, and $\chi_{-1}(x)=1-\chi(x)$; then
$$ 1=\chi_{-1}(x)+\sum_{h=0}^{\io} \chi_h(x)=\sum_{h=-1}^{\io} \chi_h(x) .
\tag 2.3$$
Calling
$$x_{n,m}(\e,\n)= |\O n|-\sqrt{\o_m^2+n\n_{n,m}}\tag 2.4$$ 
we define
$$g_{n,m,h}= \chi_{h}(x_{n,m}(\e,\n)) g_{n,m}(\e,\n)\tag 2.5$$ 
Note that if $\chi_{h}(x) \neq 0$ for $h\ge 0$ one has
$2^{-h-1}\g\le|x|\le 2^{-h+1}\g$, while if $\chi_{-1}(x)\neq 0$
one has $|x| \ge \g$. Therefore   $g_{n,m,h}(\e,\n)=0$  whenever 
$2^{-h-1}\g\le|x_{n,m}(\e,\n)|\le 2^{-h+1}\g$ is not verified.
Moreover if $g_{n,m,h}(\e,\n)\neq 0$ and  $g_{n,m,h'}(\e,\n)\neq 0$ then necessarily $|h-h'|\le 1$.
Inserting  2.3 in  2.2 we get
$$u_{n,m}^{(k)} =\sum_{h}g_{n,m,h}F_{n,m}^{(k)}+ n\sum_{h=-1}^\io
\sum_{r=2}^{k-1} l_{n,m}^{(r)} g_{n,m,h}u^{(k-r)}_{n,m}\equiv\sum_{h=-1}^\io u_{n,m,h}^{(k)}
\tag 2.6$$
\smallskip 
\noindent{\bf 2.2)}  2.6 can be applied recursively until we obtain $u_{n,m}^\ka$ 
as a (formal) polynomial in the variables
$g_{n,m,h}$, $q$  and $l^{(r)}_{n,m} $ with $r<k$. It turns out that 
$u_{n,m}^\ka$ can be written as sum over {\it trees} (see Lemma 3 below) defined
in the following way.

\0A (connected) graph $\GG$ is a collection of points (vertices)
and lines connecting all of them. The points of a graph
are most commonly known as graph vertices, but may also be
called nodes or points. Similarly, the lines connecting
the vertices of a graph are most commonly
known as graph edges, but may also be called branches or
simply lines, as we shall do. We denote with
$V(\GG)$ and $L(\GG)$ the set of vertices (also called nodes) and the set of lines,
respectively. A path between two vertices is a subset of $L(\GG)$
connecting the two vertices. A graph is planar if it can be drawn
in a plane without graph lines crossing.
\smallskip

\0{\bf Definition 3.}
{\it A {\rm tree} is a planar graph $\GG$ containing no closed loops
(cycles).
One can consider a tree $\GG$ with a single special vertex $v_{0}$:
this introduces a natural partial ordering on the set
of lines and vertices, and one can imagine that each line
carries an arrow pointing toward the vertex $v_{0}$.
We can add an extra (oriented) line $\ell_{0}$
exiting the special vertex $v_{0}$; the added line will
be called the {\rm root line}. In this way we obtain
a {\rm rooted tree} $\th$ defined by $V(\th)=V(\GG)$
and $L(\th)=L(\GG)\cup\ell_{0}$. A {\rm labeled tree} is a
rooted tree $\th$ together with a label function defined on
the sets $L(\th)$ and $V(\th)$.}

\smallskip

We shall call {\it equivalent} two rooted trees which can be transformed
into each other by continuously deforming the lines in the plane
in such a way that the latter do not cross each other
(i.e. without destroying the graph structure).
We can extend the notion of equivalence also to labeled trees,
simply by considering equivalent two labeled trees if they
can be transformed into each other in such a way that also
the labels match.

Given two {\it nodes} (sometimes also called {\it vertices}) 
$v,w\in V(\th)$, we say that $w \prec v$
if $v$ is on the path connecting $w$ to the root line.
We can identify a line with the nodes it connects;
given a line $\ell=(v,w)$ we say that $\ell$
enters $v$ and comes out of $w$. 

In the following we shall deal mostly with labeled trees:
for simplicity, where no confusion can arise, we shall call them just
trees. 
\smallskip

\0 We call {\it internal nodes} the vertices such that there is at least
one line entering them. We call {\it end-points} the vertices
which have no entering line. We denote with $L(\th)$, $V_0(\th)$
and $E(\th)$ the set of lines, internal nodes and end-points, respectively.
Of course $V(\th)=V_0(\th)\cup E(\th)$.
\smallskip
We call  $\Theta^{(k)}_{ n,m}$ the set of all the possible trees {\it of order $k$} 
defined according to the following rules.

\bigskip
\insertplotttt{197pt}{105pt}{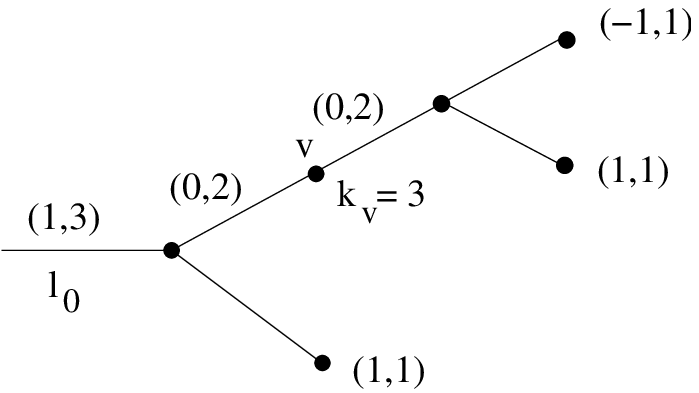}
\medskip
\line{\vtop{\line{\hskip1.5truecm\vbox{\advance\hsize by -3.1 truecm
\noindent{{\bf Fig. 1}. A tree $\th\in\Th_{3,4}^{(4)}$
}} \hfill} }}
\bigskip
\smallskip\0
(1) To each end-point $v\in E(\th)$ one associates 
the mode label $(n_{v},m_{v})$, with $m_{v}= 1$ and
$n_{v}=\pm 1$, such that
$$ \sum_{v\in E(\th)} n_{v}= n. \tag 2.7 $$
we associate to each end-node a factor $\h_v=q$ and an order $k_v=0$.

\smallskip\0

(2) To each line $\ell\in L(\th)$  one associates the mode label $(n_\ell,m_\ell)$ where one has
$$  n_{\ell} = \sum_{w \in E_{\ell}}  n_{w} 
\tag 2.8 $$
where $E_\ell$ are the endpoints of the subtree with root given by $\ell$.

(3) To each line $\ell\in L(\th)$  one associates the scale label $h_{\ell}\in \NNN \cup\{-1,0\}.$ If two lines $\ell,\ell'$ have the same mode label $(n_\ell,m_\ell)=(n_\ell',m_\ell')$ then $|h_\ell-h_{\ell'}|\le 1$. If $\ell$ exits an end-node then $h_\ell=-1$.

(4) To each node $v\in V_0(\th)$ is associated a type label $t_v=a$ or $b$;  For each node $v\in V_0(\th)$ one has $s_v=1,2$ entering lines.

If $s_v=1$ the momenta of the exiting and entering line
are necessarily the same and the type label is by definition $a$.
To $v$ is associated
an order $k_v\in [2,\io)$ and a factor $\h_{v}=n_\ell l^{(k_v)}_{n_\ell, m_{\ell},h_\ell}$ where $\ell$ is the line exiting $v$.

If $s_v=2$ then necessarily $k_v=1$. Calling $m,m_1,m_2$ 
the momenta $m_\ell$ respectively of the lines exiting and entering  $v$,  
to $v$ is associated a factor 
$\h_v=a v_{m,m_1,m_2}$ if $v$ is of type $a$ and 
$\h_v=-b\O^2v_{m,m_1,m_2}$if $v$ is of type $b$.

(5) To each line entering an $a$ node and to the root line of each tree, we associate the {\it propagator}
$$ g_{\ell} \= g_{n_{\ell} ,m_{\ell},h_{\ell}}(\e,\n) =
\left\{\eqalign{
{\chi^{(h_\ell)}(|\O n_{\ell}|-\sqrt{\o_{m_\ell}^2+n_\ell\n_{n_\ell,m_\ell}})\over -\O^{2} n_{\ell}^2 +\o_{m_\ell}^2+
n_\ell\n_{n_\ell,m_\ell}} , &
(n_{\ell},m_{\ell}) \neq (\pm 1,  1), \cr
& \cr
1  & (n_{\ell},m_{\ell}) = (\pm 1, 1). \cr}\right.
\tag 2.9
 $$

To each line entering a $b$ node we associate 

$$ g_{\ell} \=   n_{\ell}g_{n_{\ell} ,m_{\ell},h_{\ell}}(\e,\n)=
\left\{\eqalign{
{ n_{\ell}\chi^{(h_\ell)}(|\O n_{\ell}|-\sqrt{\o_{m_\ell}^2+n_\ell\n_{n_\ell,m_\ell}}) \over -\O^{2} n_{\ell}^2 +\o_{m_{\ell}}^2+n_\ell\n_{n_{\ell} ,m_{\ell}}} &
(n_{\ell},m_{\ell}) \neq (\pm 1,  1), \cr
& \cr
1   &(n_{\ell},m_{\ell}) = (\pm 1, 1). \cr}\right.
\tag 2.10 $$
Only the lines coming out from
the end-points can have momentum $(n_{\ell},m_{\ell})=(\pm 1, 1)$.

(6) Finally we define the order of a tree as:
$$k(\th)= \sum_{v\in V(T)}k_v . \tag 2.11$$
 Note that $|n_\ell|< k(\th)-\sum_{v\atop s_v=1}k_v$.
\smallskip
By the support properties of $\chi_h$ 
and bounding the denominator of $g_\ell$ with $C |n_\ell| 2^{-h_\ell}$, we get
%
$$|g^{\ell}|\le C 2^{h_{\ell}+1}\tag 2.12$$
The divisors can be small only if $n_\ell\simeq m_\ell^2$,
as explained by the following Lemma.
\smallskip
\noindent{\bf Lemma 2. }{\it If $g^{\ell}\not= 0$ and $h_\ell\ge 0$ then
$$|\o_1 |n_\ell|-m_\ell^2|\le 1+\e_0 |n_\ell|\tag 2.13$$}
\smallskip
{\it Proof.} Equation
 2.13 is equivalent to $ (\o_1-\e_0)
|n_\ell|-1\le m_\ell^2\le (\o_1+\e_0)|n_\ell|+1$;
we claim that if  $m_\ell^2> (\o_1+\e_0)|n_\ell|+1 $ or $m_\ell^2< 
(\o_1-\e_0)|n_\ell|-1$ 
then  $n_\ell\not =[\O^{-1} m^2_\ell]$ ($[...]$ denotes the closest integer); 
in fact if 
$n_\ell=\O^{-1} m_\ell^2+ x$ with $|x|\leq {1\over 2}$ then as 
$\o_1-\e_0<\O<\o_1+\e_0$ we have that:
$$ (\o_1-\e_0)(|n_\ell|-{1\over 2})\leq m_\ell^2 \leq (\o_1+\e_0)(|n_\ell|+{1\over 2})\tag 2.14. $$
as $\o_1+\e_0\le 2$. 

By contradiction assume that  2.13 is not true; then $n_\ell\not =[\O^{-1} m^2_\ell]$; then
we can write
$n= \O^{-1} m_\ell^2+ k+ x$ with $|x|\leq {1\over 2}$, $|k|\ge 1$ so that
%
$$|\O n_\ell-\sqrt{\o_{m_\ell}^2+n_\ell\n_{n_\ell,m_\ell}}|
\ge |\O n_\ell-m_\ell^2|-|m_\ell^2-\sqrt{\o_{m_\ell}^2+n_{\ell}
\n_{n_\ell,m_\ell}}|\ge $$ 
$$|\O n_\ell-m_\ell^2|- {\m_0+|n_\ell\n|_\io\over m_\ell^2}\ge  |k|-{1\over 2}- m_\ell^{-2}(\m_0+(\O^{-1} m_\ell^2+ 
k+{1\over 2})\e_0)\ge {1\over 8}>\g\tag 2.15$$ 
in contradiction with  $g^{\ell}\not= 0$ and $h_\ell\ge 0$.
\qed\
\smallskip
The coefficients $u^{(k)}_{ n,m}$ can be represented as sum over the trees
defined above; this is in fact the content of the following Lemma. 
\smallskip
\noindent{\bf Lemma 3.} {\it $u^{(k)}_{ n,m}$ solving  2.6 can be written as
$$ u^{(k)}_{ n,m} = 
\sum_{\th\in\Theta^{(k)}_{ n,m}} \Val(\th) ,
\tag 2.16$$
where
$$ \Val(\th) = 
\Big( \prod_{\ell \in L(\th)} g_{\ell} \Big)
\Big( \prod_{v \in V(\th)} \h_{{v}} \Big)\;.
\tag 2.17 $$}
\smallskip
{\it Proof.} The proof is done by induction on $k$.  
If $k=1$ it holds by  2.6, recalling that 
$u^{(0)}_{n,m}=q \d_{n,\pm 1}\d_{m, 1}$ (see Fig.2)
$$u_{n,m}^{(1)}=\sum_{h=-1}^\io g_{n,m,h}
\sum_{n_1=\pm 1} v_{m,1,1} (a-b\O^2 n_1 (n-n_1))u^{(0)}_{n_1,1}
u^{(0)}_{n-n_1,1}=\sum_{\th\in \Th^{(1)}_{n,m}}\Val(\th)\;.\tag 2.18$$
%

\bigskip
\insertplotttt{200pt}{60pt}{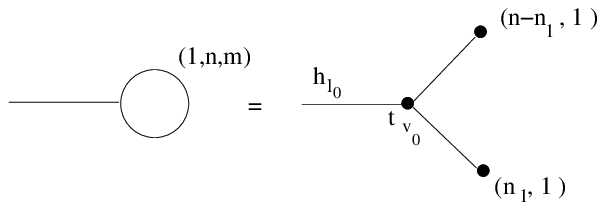}
\line{\vtop{\line{\hskip1.5truecm\vbox{\advance\hsize by -3.1 truecm
\noindent{{\bf Fig. 2}.Graphical representation of  2.18 for $k=1$; the sum over $n_1,h_{\ell_0}, t_{v_0}$ is understood. 

}} \hfill} }}
\bigskip
\smallskip\0
>From  2.6,  2.2 and the inductive hypothesis we have
that $u^{(k)}_{n,m}$ is given by
$$u^{(k)}_{n,m}=\sum_{h=-1}^\io g_{n,m,h}\{\sum_{n_1,m_1,m_2,k_1}
v_{m,m_1,m_2} $$
$$(a- b\O^2 n_1 (n-n_1))
\sum_{\th_1\in \Th^{(k_1)}_{n_1,m_1}}\Val(\th_1)
\sum_{\th_2\in \Th^{(k-k_1)}_{n-n_1,m_2}}\Val(\th_2)
+ n \sum_{r=2}^{k-1} l_{n,m}^{(r)}
\sum_{\th_3\in \Th^{(k-r)}_{n,m}}\Val(\th_3)\}\tag 2.19$$
which can be expressed graphically from Fig.3. 
\medskip
\insertplotttt{200pt}{60pt}{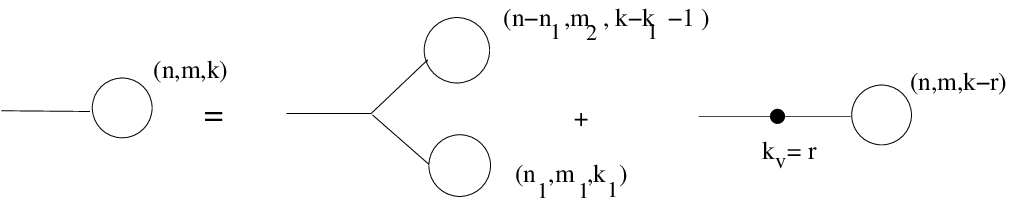}
\medskip
\line{\vtop{\line{\hskip1.5truecm\vbox{\advance\hsize by -3.1 truecm
\noindent{{\bf Fig. 3}.Graphical representation of  2.19; the sum over $n_1,m_1,k_1,m_2,r$ is understood.
}} \hfill} }}
\bigskip

Given a tree $\th\in \Th^{(k)}_{n,m}$ such that $s_{v_0}=2$, $h_{l_0}=h$,
let  $\th_1\in \Th_{n_1,m_1}^{(k_1)}$, $\th_2\in \Th_{n-n_1,m_2}^{(k-k_1)}$ 
be the subtrees whoose root lines enter in $v_{0}$; if $v_0$ is of type $a$
by   2.17 one has that:
$$\Val(\th)= a g_{n,m,h} v(m,m_1,m_2)\Val(\th_1)\Val(\th_2);\tag 2.20 $$ 
if $v_0$ is of type $b$ let $n_1$ be the momentum of $\th_1$.
By our definitions we have that:
$$\Val(\th)= -b\O^2 n_1(n-n_1) g_{n,m,h} v(m,m_1,m_2)\Val(\th_1)\Val(\th_2);\tag 2.21. $$ (recall that the root line of a tree is always an a-line.)
Finally  given  a tree $\th\in \Th^{(k)}_{n,m}$ such that $s_{v_0}=1$, 
$k_{v_0}=r$, $h_{l_0}=h$ let $\th_3$ be the subtree whoose root line enters $v_{0}$, by  2.17 one has that:
$$ \Val(\th)=n g_{n,m,h} l_{n,m}^{(r)}\Val(\th_3)\;.\tag 2.22$$ 
Hence inserting  2.20,  2.22 in  2.19 we get  2.16.
\qed
\smallskip
\bigskip
\noindent{\bf 3. Choice of the parameters $l_{n,m}$.}
\smallskip 
\noindent{\bf 3.1)} 
In the preceding section we have found a power series expansion for $u_{n,m}$
solving  1.11 and parametrized by $l_{n,m}$.
However for generic values of $l_{n,m}$ 
such expansion is not convergent, as one can easily identify contributions
at order $k$ which are $O(k!^\a)$, for a suitable constant $\a$.
In this section we show that it is possible to choose the parameters
$l_{n,m}$ in a proper way to cancel such ``dangerous'' contributions;
in order to do this we have to identify the dangerous contributions
and this will be done through the notion of {\it clusters} and 
{\it resonances}.
\smallskip
\noindent{\bf Definition 4.}{\it Given a tree $\th\in\Th^{(k)}_{n,m}$ a cluster $T$ is a connected 
set of nodes which are linked by a continuous 
path of lines with the same scale label $h_{T}$ or
a lower one and which are maximal; we shall say that the cluster
has scale $h_{T}$. We shall denote by $V(T)$ and $E(T)$
the set of nodes  and the set of end-points, respectively, which are
contained inside the cluster $T$, and with $L(T)$
the set of lines connecting them.}
\smallskip
Therefore an inclusion relation is established between clusters,
in such a way that the innermost clusters are the clusters
with lowest scale, and so on.
Each cluster $T$ has an arbitrary number of lines entering it (incoming
lines), but only one or zero line coming from it (outcoming or {\it root} line);
we shall denote the latter (when it exists) with $\ell_{T}^{1}$, and
we shall denote by $h_{T}^{(e)}$
the scale of the outcoming external line of $T$.

\smallskip
\noindent{\bf Definition 5.} {\it A cluster $T$ with $|V(T)|>1$, with
only one incoming line $\ell_{T}$
such that one has
$$  n_{\ell_{T}^{1}} =  n_{\ell_{T}} \hbox{  and  }
m_{\ell_{T}^{1}}=  m_{\ell_{T}}
\tag 3.1 $$
will be called {\it resonance} of scale $h$.
In such a case we shall call a {\it resonant line}
the root line $\ell_{T}^{1}$. 
}
\smallskip 
The propagators on the path between the external lines
of $T$ have the form, $\a_\ell=(0,1)$
$${n_\ell^{\a_\ell} \chi^{(h_\ell)}(|\O n^0_\ell+x|-\sqrt{\o_{m_\ell}^2+n_\ell\n_{n_\ell,m_\ell}}) 
\over -(\O n^0_\ell+x)^2+\o_{m_\ell}^2+n_\ell\n_{n_\ell,m_\ell}}\Big|_{x=\O n_{\ell_T}}\tag 3.2$$
and we can consider the value of $T$ as a function of $m,n,x=\O n_{\ell_T}$.
The contribution of a resonance $T$ of a tree $\th$ is given by, calling
$(n_{\ell_T},m_{\ell_T})=(n,m)$:
$$\VV_{T}^h(\O n,m,n)=
\Big( \prod_{\ell \in T} g_{\ell} \Big)
\Big( \prod_{v \in V(T)} \h_{v} \Big) .
\tag 3.3 $$
with $h=h^{(e)}_T$. 
\smallskip
We define the localization operation acting on the resonances $T$
in the following way; if $|\o_1|n|-m^2|\le 1+\e_0|n_0|$
and $(n_\ell,m_\ell)\neq (n,m)$ for all $\ell\in T$ then
$$\LL\VV_T^h(\O n,m,n)=\VV_T^h(\sign(n)\sqrt{(\o_m^2+n\n_{n,m}}),m,n)\tag 3.4$$
and $\LL=0$ otherwise. We split each resonance as
$$\VV_T^h(\O n,m,n)=\LL\VV_T^h(\O n,m,n)+\RR \VV_T^h(\O n,m,n)\tag 3.5$$
where $\RR=1-\LL$;
we call $\LL\VV_T^h(\O n,m)$ {\it local resonances}. The action of $\LL$
is then to replace, in the path connecting the external lines of $T$, the variable

$x$ with 
$$\bar\o_{m,n}=\sign(n)\sqrt{(\o_m^2+n\n_{n,m}})\tag 3.6$$
\smallskip
\noindent{\bf Definition 6.} {\it The trees $\th_T\in \R^\ka_{h,n,m}$ 
are defined as the trees $\th\in \Th^\ka_{h,n,m}$
with the following modifications: a)
there is a single end node, called $e$,
such that 
$(n_e,m_e)= (n,m)\not=(\pm 1,1)$; to $e$ is associated
$\h_e= 1/m_e^3$. If 
$\ell_e$ be the line exiting from $e$, $\ell_e$ has associated $g_{\ell_e}=1$ if it enters an $a$
 node and $g_{\ell_e}= n_e$ if it enters a $b$ node;
  
\0 b)the root line $l_0$ has $(n_{\ell_0},m_{\ell_0})= (n,m) $ and $g_{\ell_0}=1$;

\0 c) for all lines $\ell\in \th$: $\max_{\ell\in L(\th)\setminus \{\ell_0,\ell_e)\}}$ $(h_\ell)=h$.}

\smallskip
The definition of value of such tree is identical to the one given in  2.17

\bigskip
\insertplotttt{297pt}{105pt}{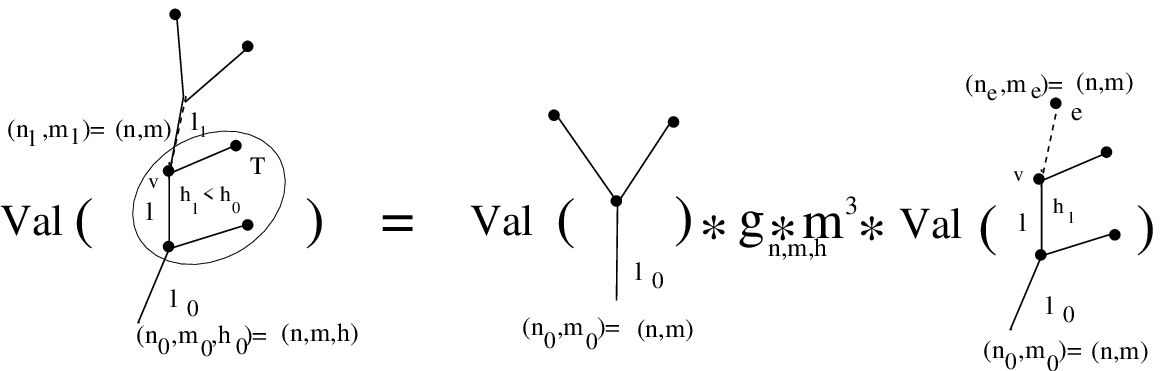}
\medskip
\line{\vtop{\line{\hskip1.5truecm\vbox{\advance\hsize by -3.1 truecm
\noindent{{\bf Fig. 4}.
We 
associate to the resonance $T$ (enclosed in an ellipse) the tree $\th_T\in \R$, and vice-versa. 
}} \hfill} }}
\medskip

Given a resonance $T$, there exists a 
unique $\th_T\in\R^\ka_{h,n,m}$ such that (see Fig 3.1)
$$\VV_{T}(\O n,m,n)=m^3  \Val(\th_T)\tag 3.7$$ 
where $\th_T\in\R_{n,m,h}$ if the external line enters an $a$ node and 
$n \VV_{T}(\O n,m,n)=m^3  \Val(\th_T) $ if the external line enters an $b$ node.

\smallskip 
\noindent{\bf 3.2}
With a suitable choice of the parameters $l_{n,m}$ the
functions $u_{n,m}^{(k)}$ can be rewritten as sum over 
``renormalized'' trees defined below.
\smallskip
\noindent{\bf Definition 7.}
{\it We define the set of {\it renormalized tree}
$\Th_{R,n,m}^{(k)}$ defined as the trees in 
$\Th_{R,n,m}^{(k)}$ defined in \S 2 with the following differences:
a)to each resonance $T$ we apply the $\RR$ operation; b)the nodes with 
$s_v=1$ have associated 
$\h_v= n_{\ell}l^{(k_v)}_{n_\ell,m_\ell,h_\ell}$ 
where $\ell$ is the line entering $v$.
In the same way we define 
$\R_{R,h,n,m}^{(k)}$. We call resonant lines the lines coming out a resonance or a node with 
$s_v=1$.}
\smallskip
It holds the following result.
\smallskip
\noindent{\bf Lemma 4.}{\it For all $k,n,m$ it holds:
$$u_{n,m}^{(k)}= \sum_{\th\in \Th_{R,n,m}^{(k)}}\Val(\th)\tag 3.8$$
with
$$n l_{n,m,h}^{(k)}=-m^3 \sum_{h_1\ge h}
\sum_{\th\in \R_{R,n,m,h_1}^{(k)}
}\LL\Val(\th)\tag 3.9$$ 
provided that we choose $l_{n,m}^{(k)}=
l_{n,m,-1}^{(k)}$  in  2.6.
}
\smallskip
\medskip
\insertplotttt{237pt}{75pt}{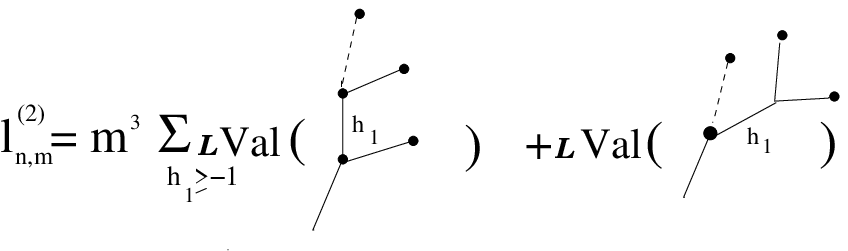}
\line{\vtop{\line{\hskip1.5truecm\vbox{\advance\hsize by -3.1 truecm
\noindent{{\bf Fig. 5}.
The counterterm $l_{n,m}^{(2)}$.
}} \hfill} }}
\bigskip

{\it Proof.} First note that by definition
$l_{n,m,h}= 0$ if $|\o_1 n- m^2|\ge 1+\e_0|n| $.
We proceed by induction. For $k=1,2$ 3.8 surely 
holds as $\Th_{R,n,m}^{(1,2)}\equiv\Th_{n,m}^{(1,2)}$.
Then we assume that 3.8 holds for all $r<k$; by  2.6
$$u_{n,m,h}^{(k)}=  g_{n,m,h} 
(F_{n,m}^\ka +n\sum_{r=2}^{k-1} l^{(r)}_{n,m}u^{(k-r)}_{n,m})\tag 3.10$$ 

$F^\ka_{n,m}$ is a function of $u_{n'.m',h'}^{(r')}$ with $r'<k$ which,
by the inductive hypothesis, are written as sum over trees 
in $\Th_{R,n',m'}$.
$g_{n,m,h} F_{n,m}^\ka$ is given by sum over 
$\th\in \Th^{(k)}_{n,m}$ with $s_{v_0}=2$, and the root line can 
be resonant or not. If $l_0$ is non-resonant then $\th\in \Th_{R,n,m}^\ka$. If $l_0$ is resonant
we split the biggest resonance in the form  3.5; if $\LL=0$
necessarily there is an inner resonance (whose resonant line is the root
line) and again we apply  3.5 and surely $\LL\not=0$. We split
$g_{n,m,h} F_{n,m}^\ka$ as sum of two terms; one, which we denote by 
$G_{n,m}^{(k)}$, which is the sum over all trees 
belonging to $\Th_{R,n,m}$ with $s_{v_0}=2$ and the second
which is sum of trees with value
$$\Val(\th)=g_{n,m,h_{\ell_0}}[\LL\Val(\th_T)]\Val(\th_1)\tag 3.11$$
with $\th_T\in \R^{(r)}_{R,h_1,n,m}$ and 
$\th_1\in \Th^{(k-r)}_{R,n,m}$. We get
$$F_{n,m}^\ka = m^3  g_{n,m,h}
\sum_{r=2}^{k-1}u_{n,m}^{k-r}\sum_{h_1<h} 
(\LL\sum_{\th\in \R_{R,n,m,h_1}^{(r)}}\Val(\th))\, + G_{n,m}^\ka \tag 3.12$$ 
which inserted in  3.10 and using  3.9 gives
$$u_{n,m,h}^{(k)}= g_{n,m,h} m^3
\sum_{r=2}^{k-1}u^{(k-r)}_{n,m}(\sum_{h_1<h} 
\sum_{\th\in \R_{R,n,m,h_1}^{(r)}}\LL\Val(\th))\, + G_{n,m}^\ka \tag 3.13$$
$$-m^3\sum_{r=2}^{k-1}u^{(k-r)}_{n,m}
(\sum_{h_1\ge -1}\sum_{\th\in 
\R_{R,n,m,h_1}^{(r)}}\LL\Val(\th))= g_{n,m,h}G_{n,m}^\ka
+g_{n,m,h}\sum_{r=2}^{k-1}u^{(k-r)}_{n,m} l_{n,m,h}
$$
and by definition $G_{n,m}^\ka$ is a sum over all $\th \in\Th^{(k)}_{R,n,m}$.
\qed
\bigskip
\insertplotttt{350pt}{100pt}{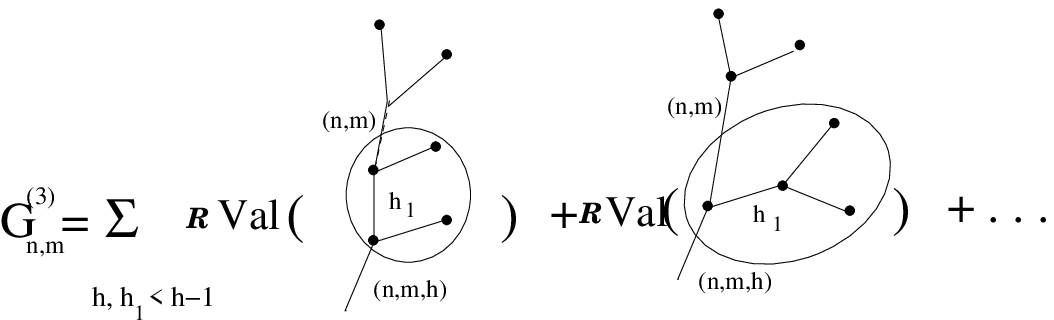}
\medskip
\line{\vtop{\line{\hskip1.5truecm\vbox{\advance\hsize by -3.1 truecm
\noindent{{\bf Fig. 6}.The term  $G_{n,m}^{(3)}$, the dots represent sums over trees
 with $s_v=2$ and  non resonant root line.
}} \hfill} }}
\bigskip


\smallskip
%
%
\smallskip

\bigskip
\noindent{\bf 4. Bruno Lemmas and bounds for the expansion}
\smallskip 
\noindent{\bf 4.1)} In the previous section we have shown that, with a suitable choice of the parameters 
$l_{n,m}$, we can express $u_{n,m}$ as sum over trees belonging to $\Th^k_{R,n,m}$;
we show in this section that such expansion is indeed convergent
if $\h$ is small enough and $\e,\n\in D(\g)$ (see Definition 1).

Given a tree $\th\in \Th_R$, we call  $S(\th,\g)$ the set of  $(\e,\n)\in \D$ such that:
 for all $\ell\in L(\th)$: $$2^{-h_\ell-2}\g< |\O |n_\ell|- \sqrt{\o_{m_\ell}^2+n_\ell \n_{n_\ell,m_\ell}}|<
2^{-h_\ell+2}\g.\tag 4.1$$

In other words we can have $\Val(\th)(\e,\n)\neq 0$   only if $(\e,\n)\in S(\th,\g)$.

 We call  $D(\th,\g)\subset \D$  the set of $(\e,\n)$ such that, if 
$\a_1=\pm$, $\a_2=\pm$:
$$|x_{n_\ell,m_\ell}|=\left| \O |n_\ell|- \sqrt{\o_{m_\ell}^2+n_\ell\n_{n_\ell,m_\ell}} \right|
\ge \g |n_\ell|^{-\t}\tag 4.2
$$
$$|x^{\a_1,\a_2} _{n_{l_1},m_{l_1},n_{l_2},m_{l_2}}|=$$
$$\left| \O (n_{\ell_1}-n_{\ell_2})+ \a_1 
\sqrt{\o_{\ell_1}^2+n_{\ell_1} \n_{n_{\ell_1},m_{\ell_1}}}+
\a_2 \sqrt{\o_{\ell_2}^2+n_{\ell_2} \n_{n_{\ell_2},m_{\ell_2}}}
\right| \ge \g |n_{\ell_1}-n_{\ell_2}|^{-\t}$$
$$\forall   |\o_1 n_{\ell_i}- m^2_{\ell_i}| <1 +\e_0 |n_{\ell_i}|\tag 4.3$$ 
for all lines $\ell_1,\ell_2\in L(\th)$ such that $n_{\ell_1}\neq n_{\ell_2}$
This means $D(\th,\g)$ is the set of $(\e,\n)$ verifying the Melnikov conditions in $\th$.

Calling $L_0(\th),V_0(\th)$
the set of lines, node and end-points 
not contained in any resonance, and $S_{0}(\th)$
the {\it maximal resonance}, \ie the resonances
which are not contained in any other resonance,
we can write $\Val(\th)$ with $\th\in  \Th_{R,n,m}$ as 
$$ \Val(\th) = 
\Big( \prod_{\ell \in L_0(\th)} g_{\ell}^{(h_{\ell})} \Big)
\Big( \prod_{v \in V_0(\th)} \h_{{v}} \Big) 
\Big(\prod_{T\in S_0(\th)} \RR \VV^{h_T^{e}}_T(\O n_{\ell_T},m_{\ell_T}, n_{\ell_T})
\Big) ,
\tag 4.4 $$
and by definition
$$ \RR \VV^{h_T^{e}}_T(\O n_{\ell_T},m_{\ell_T},n_{\ell_T})=$$
$$ \VV^{h_T^{e}}_T(\O n_{\ell_T},m_{\ell_T},n_{\ell_T})-
\VV_T^{h_T^{e}}({\text  sign}(n_{\ell_T})\,\sqrt{(\o_{m_{\ell_T}}^2+n_{\ell_T}\n_{n_{\ell_T},m_{\ell_T}}})
,m_{\ell_T}, n_{\ell_T}) ,
\tag 4.5 $$
and $\VV^{h_T^{e}}_T(\O n_{\ell_T},m_{\ell_T},n_{\ell_T})$
is given by
$$\VV^{h_T^{e}}_T(\O n_{\ell_T},m_{\ell_T},n_{\ell_T})=
\Big( \prod_{\ell \in L_0(T)} g^{(h_{\ell})}_{\ell} \Big)
\Big( \prod_{v \in V_0(T)} \h_{v} \Big)
\Big(\prod_{T'\in S_0(T)} \RR \VV^{h_{T'}^{e}}_{T'}
(\o n_{\ell_{T'}},m_{\ell_{T'}}) \Big) .
\tag(4.6) $$
In order to bound $\Val(\th)$  4.4 we will use the following result.
\smallskip
\0\noindent{\bf Lemma 5 (Bruno Lemma).} {\it
Given  tree $\th \in \Th^{(k)}_{R,n,m}$, 
we have that  $D(\th,\g)\cap S(\th,\g)\neq 0$ if and only if the scales 
$h_\ell$ of $\th$ respect
$$ N_h(\th) \le  \max\{ 0 , 2 K(\th) 2^{(2-h)/\t} - 1 \}+S_h(\th)+M_h(\th) , 
\tag 4.7 $$
where $ N_h(\th) $ is the number of lines with scale greater or equal than 
$h$, $K(\th)\leq k(\th)$ is the number of non resonant lines,
 $S_h(\th)$ is the number of resonances $T$
in $\th$ with $h_{T}^{(e)}=h$
and $M_{h}(\th)$ is the number of vertices with $s_v=1$ in $\th$
such that the scale of the exiting line is $h$.}
\smallskip
The proof of the above Lemma is in Appendix A3. By the above lemma 
we can prove the following result.
\smallskip
\0{\bf Lemma 6.} {\it Assume that there exist a constant $C$ such that
one has $|l^{(k)}_{h,n,m}|\le q^{2k}C^{k-1} 2^{-h}$, for any $n,m$
and all $h\ge 0$. Then for all 
$(\e,\n)\in \D(\th,\g)$ it holds that, for a suitable constant $D$
$$|\Val(\th)|\le D^k q^{2k}\Big(\prod_{v \in V(T)\atop s_v=2} |\h_{v}|\Big)\tag 4.8$$}
\smallskip
{\it Proof.}
Consider a tree with fixed scales $h_{\ell}$ and momenta $n_\ell,m_\ell$.
In order to take into account the $\RR$
operation we write  4.5 as, if $\bar\o_{n,m}=\sign(n)\sqrt{\o_{m}^2+n \n_{n,m}}$
$$\align &\RR \VV^{h_T^{e}}_T(\O n_{\ell_T},m_{\ell_T},n_{\ell_T})=\\ &
\left(\O n_{\ell_T}-\bar\o_{n_{\ell_T}, m_{\ell_T}}
\right)
\int_{0}^{1} {\text d}t \dpr
\VV^{h_T^{e}}_T(\O n_{\ell_T}+t(\O n_{\ell_T}-\bar\o_{n_{\ell_T},m_{\ell_T}}) ,
m_{\ell_T},n_{\ell_T}) ,
\tag 4.9 \endalign $$
where $\dpr$ denotes the derivative with respect to the argument
$\o n_{\ell_T}+t(\o n_{\ell_T}-\tilde\o_{m_{\ell_T}})$.

By (4.6) we see that the derivatives can be applied
either on the propagators in $L_0(T)$, or 
on the $\RR \VV^{h_{T'}^{e}}_{T'}$.
In the first case there is a factor 
$2^{-h^{(e)}_T+h_{T}}$: $2^{-h^{(e)}_T}$ is obtained from
$\o n_{\ell_T}-\bar\o_{n_{\ell_T},m_{\ell_T}}$
while $\partial g^{(h_T)}$ is bounded proportionally to $2^{2 h_T}$;
in the second case note that $\dpr_t \RR \VV^{h_{T'}^{e}}_{T'}=
\dpr_t \VV^{h_{T'}^{e}}_{T'}$
as $\LL\VV_{T'}^{h^{(e)}_{T'}}$ is independent of
$t$; if the derivative acts on the propagator
of a line $\ell\in L(T)$, we get a gain factor
$$ 2^{-h^{(e)}_T+h_{T'}} \le 2^{-h^{(e)}_T+h_{T}}
2^{-h^{(e)}_{T'}+h_{T'}} ,
\tag 4.10 $$
as $h^{(e)}_{T'}\le h_{T}$. We can iterate
this procedure until all the $\RR$
operations are applied on propagators;
at the end (i) the propagators are derived at most one time;
(ii) the number of terms so generated is $\le k$;
(iii) to each resonance $T$ a factor
$2^{-h^{(e)}_T+h_{T}}$ is associated.

Assuming that $|l^{(k)}_{h,n,m}|\le q^{2k}C^{k-1} 2^{-h}$ with $\g C>1$ and recalling definition  2.11,
for any $\th$ one obtains:
$$ \eqalign{
& \left| \Val(\th) \right| \le
C^{k} q^{2k}
\Big( \prod_{h=0}^{\io} \exp \Big[ h \log 2 \Big( 4 k
2^{-(h-2)/\t} + S_{h}(\th) +
M_{h}(\th) \Big) \Big] \Big) \cr
& \qquad \qquad
\Big( \prod_{T\atop h^{(e)}_T \ge 3 } 
2^{-h^{(e)}_T+h_T} \Big)
\Big( \prod_{h=0}^{\io} 2^{-h M_{h}(\th)} \Big)\Big(\prod_{v\in V(\th)\atop s_v=2} |\h_v|\Big)
\cr}
\tag 4.11 $$
where the first factor
is a bound for $\prod_{h} 2^{h N_h(\th)}$;
moreover $ \prod_{h=0}^\io 2^{-h M_\n^h(\th)}$
takes into account the factors $2^{-h}$ arising from
the running coupling constants $l^{(k)}_{h,n,m}$
and the action of $\RR$ produces, as discussed above, 
the factor $\prod_{T} 2^{-h^{(e)}_T+h_T}$.
Note that
$$\big(\prod_{h=0}^{\io} 2^{h S_{h}(\th)} \Big)
\Big( \prod_{T} 2^{-h^{(e)}_T} \Big) = 1\tag 4.12$$
Moreover it holds that
$$\prod_{T\atop h^{(e)}_T \ge 3 } 
2^{h_T}\le \prod_{h=0}^\io 2^{h 4  k
2^{-(h-2)/\t}}\tag 4.13$$
as for any derivative produced by the $\RR$ operation 
and acting on a propagator at scale $h$ there is surely a non resonant
propagator at the same scale (otherwise
the maximal clusters contained in a resonance are all resonances and 
$\RR=1$). Then we can write  4.11 as
$$\left| \Val(\th) \right| \le(q^2 D_1)^{k} 2^{4k \sum_{h=0}^\io h
2^{-(h-2)/\t}}\Big(\prod_{v\in V(\th)\atop s_v=2} |\h_v|\Big)\tag 4.14$$
from which  4.8 immediately follows.\qed
\smallskip
In order to bound the factors $|\h_v|$ we will use the following result proven in the Appendix A4.
\smallskip
\0{\bf Lemma 7.}{\it For all trees $\th\in \Th_{R,n,m}^\ka,\R^\ka_{R,h,n,m}$ with $s_v=2$ for all $v$  one has that
$$\sum_{\{m_\ell\}}\prod_{v \in V(\th)} |\h_{v}|\le {\bar C_2^k\over m^3}\tag 4.15$$ where  $m$ is the
 momentum associated to the root line, $\sum_{\{m_\ell\}}$ is the sum over the values of the momentum $m_\ell$  and $C_2$ depends only on $a,b$.}
\smallskip

\smallskip
Finally we have to prove that $l^{(k)}_{n,m,h}\le C^{k-1} 2^{-h} q^{2k}$.

Given a tree $\th\in \R_{R,n,m}$, we call 
$\tilde S(\th,\g)$ set of  $(\e,\n)\in \D$ such that 4.1 holds for all $l\in L(\th)$ not on 
the path between $e$ and $v_0$ and:
$$2^{-h_\ell-1}\g< ||\O n^0_\ell-\bar \o_{n,m}|-\sqrt{\o_{m_\ell}^2
+n_\ell\n_{n_\ell,m_\ell}}|<
2^{-h_\ell+1}\tag 4.16$$ 

holds for $\ell\neq \ell_e$ on the path between $e$ and $v_0$, namely $\LL \Val(\th)=0$ outside $\tilde S(\th,\g)$.
Finally let $\tilde D(\th,\g)\subset \D$ be the set of couples $(\e,\n)$ such that  4.2 holds for all $\ell$ not in the path  connecting $e$ to $\ell_0$, and  4.3 holds for all $\ell_1,\ell_2\in L(\th)$ such that $n_{\ell_1}\neq n_{\ell_2}$ and moreover either both $\ell_1,\ell_2$ are on the path connecting $e$ to $\ell_0$ or they both arent on such path.

First of all, the following generalization of Lemma 5 holds.
\smallskip
\0{\bf Lemma 8.} {\it
Given  tree $\th \in \RR^{(k)}_{R,h,n,m}$, 
we have that $\tilde D(\th,\g)\cap \tilde S(\th,\g)\neq 0$ if and only if the 
number of lines in $\th$ with scales $h_\ell$ 
verifies
$$ N_h(\th) \le  2 (K(\th)-1) 2^{(2-h)/\t} +S_h(\th)+M_h(\th) , 
\tag 4.17 $$
.}
\smallskip
It is an immediate consequence of the previous Lemma the following result.
\smallskip
\0{\bf Lemma 9.}{\it Given a tree $\th \in \RR^{(k)}_{R,h,n,m}$, and supposing that $l_{n,m,h'}^{(r)}\leq C^{r-1} q^{2r} 2^{-h'} $ for all $r<k$ then for   
$(\e,\n)\in \tilde D(\th,\g)$ it holds that
$$|\LL \Val(\th)|\le |n| C^{k-1}q^{2k} 2^{-h}\prod_{v\in V(\th)\atop s_v=2} |\h_{v}|C_2^{-k_1}\tag 4.18$$
where $k_1$ is the number of 
lines exiting a node with $s_v=2$. }
\smallskip
{\it Proof.}  The proof is essentially identical to the one of Lemma 6; the factor $n$ comes from the definition
 of $\Val(\th)$ in the case when the external line enters a $b$ node. 
To extract the factor $2^{-h}$ we recall that 
there is  
at least a non resonant 
line $\ell\neq \ell_0$ on scale $h_{\ell}=h, h-1$ 
which does not exit a node with $s_v=1$.  
By Lemma 8  we have that 
$k_1-1\ge K(\th)-1> 2^{h-1\over \t} $,   
so that $2^{k_1-1} 2^{-h}>1$. Then    
$$|\LL \Val(\th)|\le q^{2k} 2^{-h}(2 DC_2^2)^{k_1-1} C^{k-k_1} \le 2^{-h} C^{k-1}\prod_{v\in V(\th)\atop s_v=2}
 |\h_{v}|C_2^{-k_1}\tag 4.19$$ provided that $2 C\g>1$ and we choose $C=2DC_2^2$. 
Finally the factor $D$ as in Lemma 6 is of the form 
$D=  D_1 \g^{-1}$ with $D_1>1$ a pure ($\e$ and $\g$ independent) 
constant
\qed.
\smallskip
As aconsequence of  3.9 and Lemma 7 and 9 
it follows that $|l^{(k)}_{h,n,m}|\le q^{2k}C^{k-1} 2^{-h}$. 
%
%
%
\smallskip
\noindent{\bf Lemma 10.} {\it  For $\h_0$ small enough, the  following bounds hold 
for all $(\e,\n)\in D(\g)$:
$$|l_{n,m,h}|< C_1 \h 2^{-h}\,,\quad |l_{n,m}|< C_1 \h\,
\quad |u_{n,m}|<  C_0 |\h|{e^{-\s|n|}\over m^7}\,, \;(n,m)\neq (\pm 1,1), \tag 4.20$$
where $(n',m')$ are such that $|\o_1|n'|-(m')^2|\le 1+\e_0|n'| $ as otherwise $\n_{n',m'}\equiv 0$ by definition. } 
\smallskip
{\it Proof.} By definition   $D(\g)$ is contained in all $D(\th,\g)$ and in all $\tilde D(\th,\g)$ so that we can use  Lemma 6 and Lemma 9 to bound the values of trees.
First we fix an unlabeled tree $\th$ and sum over the values of the labels. Fixed  $(\e,\n)$ and given $(n_\ell,m_\ell)$ there are only two possible values for each $h_\ell$ such that $\Val(\th)\neq 0$. 
So we can sum up on the possible scale values obtaining a factor $2^k$. 
First we fix the tree $\th\in \Th_{R,n,m,h}$ and sum up the $m_\ell$ labels as in Lemma 7, we obtain a factor $m^{-3}$. 
Then we sum up on the possible values the momentum of lines exiting an end node, we obtain $4^k$; 
finally we bound by $\bar C^k$ the number of unlabeled trees.
The bound for $l_{n,m,h}$ is obtained by:

$$|l_{n,m,h}|\le{m^3\over n}\sum_{k=2}^\io \h^{k-1}
\sum_{h_1\ge h}\sum_{\th\in \R^\ka_{R,h_1,n,m}}|\LL\Val(\th)|\le \h C_1\sum_{h_1\ge h}2^{-h_1}. \tag 4.22$$

By Lemma 1 $u_{n,m}^\ka=0$ if $|n|>k$, so that, using Lemma 7:
$$|u_{n,m}|\le \sum_{k=1}^\io \h^k|u_{n,m}^\ka |\le \sum_{k=|n|}^\io\h^k C^{k}{1\over m^3},\tag 4.21$$
In order to get a better decay in $m$ we simply note that 
if $|n|\le {m^2\over 4}$ then 
, if $\ell_0$ is the root line, $l_0$ is surely an $a$-line, $h_{\ell_0}=-1$
and $|g_{\ell_0}|\le C m^{-4}$;
if $|n|\ge {m^2\over 4}$ of course $\h^{|n|}\le C \h^{|n|\over 2}m^{-4}$.
Then we get an extra $m^{-4}$ in  4.21 so that
the bound in  4.20 is found.
\qed

\smallskip

\bigskip
\noindent{\bf 5. Whitney extension and implicit function theorems}
\smallskip 
\noindent{\bf 5.1}
In this section we extend the function $u_{n,m}$,$l_{n,m}$, defined in  $D(\g)$ to the 
larger set $\D$.

\smallskip
\noindent{\bf Lemma 11.} {\it Given $\th\in \R^\ka_{R,h,n,m}$, 
we can  extend $\Val(\th)$ to a
function, called $\Val^E(\th)$, defined and $C^\io$ in $\D$
such that  the bounds of Lemma 10 hold for any $(\e,\n)\in \D$,
 $\LL\Val(\th)=\LL\Val^E(\th)$ for any $(\e,\n)\in\D(\th,2\g)\subset\D(\th,\g)$ and finally $\LL\Val^E(\th)=0 $ for $(\e,\n)\in \D\setminus \D(\g)$.
We then define:
$$ l_{n,m,h}^{(k)E}= \sum_{h_1\ge h} \sum_{\th\in \R^\ka_{R,h_1,n,m}}\LL 
\Val^E(\th)\tag 5.1$$
and $l^E$ is $C^1$ in $(\e,\n)\in\D$ and 
$$|\partial_\e l^E_{n,m}|< C_1 \h\,,\;|\partial_{\n_{n',m'}} l^E_{n,m}|< C_1 \h,
\quad  |\sum_{(n',m')\in \L}\partial_{\n_{n',m'}} l^E_{n,m}|< C_1 \h \tag 5.2$$
In the same way, given $\th\in \Th^\ka_{R,h,n,m}$, we define the extended 
value $\Val^E(\th)$. }
\smallskip
{\it Proof.} We prove first the statement for the
(more difficult) case $\th\in \RR^{(k)}_{R,h,n,m}$.
We use the $C^\io$ compact support function $\chi(t):\R\to \R^+$, defined in the previous section. 
Recall that $\chi(t)$
equal to $0$ if $|t|<\g$ and $1$ if $|t|\ge 2\g$, and $|\partial_t\chi|\le C$.
We proceed by induction let us suppose that we have 
proved Lemma 11 for $r<k$ and therefore defined $l_{n,m,h}^{(r)E}$ for $r<k$. 
Given a tree $\th\in\R_{R,h,n,m}^{(k)}$,  
%
$$\LL\Val^E(\th)=$$ $$\prod^*_{\ell\in L(\th)}\chi(|x_{n_\ell,m_\ell}| |n_\ell|^\t)
\prod_{\a_1,\a_2}
\prod^{**}_{\ell_1,\ell_2\in L(\th)}\chi(|x^{\a_1,\a_2}_{n_{\ell_1},m_{\ell_1},
n_{\ell_2},m_{\ell_2}}| |n_{\ell_1}-n_{\ell_2}|^\t))\LL\Val(\th)
\tag 5.3$$
where $\prod^*_{\ell\in L(\th)}$ is the product on the lines not on the path between $e$ and $v_0$ and $\prod^{**}_{\ell_1,\ell_2\in L(\th)}$ is the product on the couples   $\ell_1,\ell_2\in L(\th)$such that:
$n_{\ell_1}\neq n_{\ell_2}$ and  either both $\ell_1,\ell_2$ are on the path connecting $e$ to $v_0$ or they both are not on such path. Finally in each node $v$ with $s_v=1$ we set $\h_\ell=l_{n,m,h}^E$.

1. By definition $\Val^E(\th)=\Val(\th)$ for $(\e,\n)\in \D(\th,2\g)$
as in this set the $\chi$ in the above formula are identically equal to $1$;

2. By definition supp$( \Val^E(\th))\subset \tilde\D(\th,\g)$ as 
the $\chi$ in the above formula are identically equal to $0$ in the complementary to $\tilde\D(\th,\g)$;

Finally we define
$$l^{E\ka}_{n,m,h}(\e,\n)=\sum_{h_1=h}^\io  \sum_{\th\in \R_{n,m,h_1}}\LL\Val^E(\th)(\e,\n) 
\tag 5.4$$ 
which respects the bounds in Lemma 10. 
In order to prove 5.2
we proceed by induction. Given a tree $\th\in\R^\ka_{R,h_1,n,m}$, the derivatives act  on the nodes with $s_v=1$ which carry the factor 
$l_{n',m',h'}^{(r)}$ with $r<k$ so we can apply the inductive hypothesis.
On the lines $\ell$ not on the path  $e,v_0$
we get
$$|\partial_\e g_\ell|\le C|n_\ell|2^{2h_\ell}\,,\qquad 
|\partial_{\n_{n',m'}}g_\ell|\le C 2^{2h_\ell} \tag 5.5$$ 
and we use that $|n_\ell|\le k$.
On the lines  $\ell$ on the path $e,v_0$ the propagator is
given by $\LL g_{\ell}$, defined in  3.2 with $x$ replaced by 
$\bar\o_{n,m}$, so that
$$|\partial_\e \LL g_\ell|\le |n_\ell^{(0)}|2^{2h_\ell}\,,\qquad |\partial_{\n_{n',m'}}\LL g_\ell|= C 2^{2h_\ell};
\tag 5.6$$ 
where we have used that, by definition of $\D$, $ |\o_1n'-(m')^2|\le 1+\e_0|n'|$. Finally we consider the
 derivatives of the $\chi$ functions which  
produce in the bounds a factor $|n^0_\ell|^{\t+1}$, all this factors are bounded by $k^{\t+1}\le C^k$, 
so that the derivatives of $\Val(\th)$ respect the bounds  4.18.
As this bounds are uniform (indendent from $(n,m)$) so  that
$l^E_{n,m,h}$ is $C^1$ function of $(\e,\n)$. 

Moreover $\partial_{\n_{n',m'}}l^{\ka E}_{n,m} (\e,\n)$ is non vanishing only if
$|\o_1n'-(m')^2|\le 1+\e_0|n'|$ 
and if
$|n|-k<|n'|<|n+k|$; hence 
$$ \sum_{n',m'\in \L} |\partial_{\n_{n',m'}}l^{\ka E}_{n,m} (\e,\n)|\le C_0 k^{3/2}\max_{n',m'\in 
\L}|\partial_{\n_{n',m'}}l^{\ka E}_{n,m} (\e,\n)|\le C_1^k$$
where $\L$ was defined in Definition 2.
 
In the same way for $\th\in\Th_{R,n,m}$ 
$$\Val^E(\th)=\prod_{\ell\in L(\th)}\chi(x_{n_\ell,m_\ell} |n_\ell|^\t) 
\prod_{\a_1,\a_2}\prod_{\ell_1,\ell_2\in L(\th)}\chi(x^{\a_1,\a_2}_{n_{\ell_1},m_{\ell_1}, 
n_{\ell_2},m_{\ell_2}} |n_{\ell_1}-n_{\ell_2}|^\t)\Val^E(\th)
\tag 5.7$$
 and finally $$u^E_{n,m}= \sum_{k=1}^\io \h^k \sum_{\th\in \Th_{R,n,m}}\Val(\th).\tag  5.8$$ 
\qed
\smallskip
\smallskip
\noindent{\bf 5.2)}
{\it Proof of  Proposition 1.}
Lemma 10 and Lemma 11 imply 
that  $l_{n,m}^E$ and $u^E_{n,m}$ are $C^1$ in $(\e,\n)$ for $(\e,\n)\in \D$ and analytic in $\h, q$ for  $|q|\le q_0$ and  $\h\le \e_0$ such that $ Dq_0^2\e_0\ll 1$. 

 Inserting in the first of 1.11 the expansion for $u^E$
and $l^E$ we get the following equation for $q$ 
 $$q= \sum_{k=2}^\io \h^{k-2} \sum^*_{n_1,m_1\atop n_2,m_2}\sum_{k_1+k_2=k} 
\d_{n_1+n_2-1}(a -b \O n_1n_2) v_{1,m_1,m_2}
u_{n_1,m_1}^{(k_1)E}u_{n_2,m_2}^{(k_2)E}.
\tag 5.9 $$ 
%
Indeed the leading order of (2.5r) is
$$ q=-{1\over 2}q^3 (a-b\O^2) \sum_{m} v_{1,1,m}^2(2ag_{0,m}+(a-2b\O^2)g_{2,m} )+O(\h).\tag 5.10$$

One can easily verify that for all $(\e,\n)\in \D$ $|x_{0,m}|,|x_{2,m}|> {1\over 2}$ so that $u_{0,m}=u_{0,m}^E$. 
 So (the equation for $q$ is equivalent to:
$$q= (A+O(\e_0)) q^3+\h F(q,\e,\n,\h) \tag 5.11$$ 
We then exclude those values of $a,b$ for which $A\le 0$.
Equation 5.11 is clearly invertible near $\h=0$ if $A>0$, so that we obtain $q=q(\eta,\e,\n)$ analytic in $\h$ and $C^1$ in $(e,\n)$.
 This completes the proof  of  Proposition 1. 

Notice that if $A<0$ then we would only need to consider $\O=\sqrt{1+\m-\e}$
with as usual $\e>0$.
\qed
\smallskip
\bigskip
\noindent{\bf 6. Proof of Proposition 2}

\smallskip
\noindent{\bf 6.1)} In order to prove the first part of Proposition 2,
we consider the extended  compatibility equation 1.12:
$$\n_{n,m}= \h l^E_{n,m}(\e,\n,\h)\equiv \sum_{k=2}^\io \h^{k}l_{n,m}^\ka, \tag 6.1$$ 
where we have substituted $q=q(\e,\n,\h)$.
$l^E_{n,m}(\e,\n,\h)$ is a $C^1$ function with bounded Jacobian (see 5.2) 
so that we can solve 6.1 
by the implicit function theorem for $\h<\h_0$ small enough. 
We  obtain a function $\n(\e,\h)$ defined for $\e\in (0,\e_0)$, $|\h|\leq \h_0$ and of order $\h^2$. Moreover $\n_{n,m}(\e,\h)=0$ if $|\o_1|n|-m^2|\ge 1+\e_0|n|$.
 
One can derive expression 6.1:
$$\partial_\e \n_{n,m}(\e,\h)= 
\h (\partial_\e l^E_{n,m}+\sum_{n',m'\atop |\o_1|n'|-(m')^2|\le 1+\e_0|n'|}\partial_{\n_{n',m'}} l^E_{n,m}\partial_\e \n_{n',m'}(\e,\h)) 
$$
 so that   $\n_{n,m}(\e,\h)$ is differentiable in $\e,\h$ and respects the same bounds as $l_{n,m}$ namely:
$$|\partial_\e \n_{n,m}(\e,\h)|_\io \leq \h^2 C \,,
\qquad |\partial_\h \n_{n,m}(\e,\h)|_\io \leq \h C$$
Finally we set $\h=\sqrt{\e}$ and obtain the desired bounds.\qed
\smallskip
\noindent{\bf 6.2)}
We have now to bound the measure of $\CC(\g)$.
%

We define $\II_1$ the set of $\e\in (0,\e_0)$ verifying 
for any  $(n,m)$, 
$$|\O n -\sqrt{\o_m^2+n\n_{n,m}(\e)}|\le{C_0\over |n|^{\t_0}}\tag 6.2$$
with $C_0=2\g$.
When l.l5 is satisfied by Lemma 2 there exists two constants such that
$$c_1 \sqrt{n}\le m\le c_2 \sqrt{n}\tag 6.3$$
Moreover one must have
(by using also 1.4)
$$ \eqalign{
2 C_{0} | n|^{-\t_{0}} & \le | \o_{1} n - \o_{m} | \cr
& \le | \o_{1} n + \e n -
\sqrt{\o_{m} + n\n_{n,m}(\e)}| \cr
& \qquad \qquad +
| \e n|+ |\sqrt{\o_m^2+n \n_{m}(\e)}-\o_m| \cr
& \le C_{0} | n|^{-\t} + C\e_{0}| n| , \cr}
\tag 6.4 $$
which implies, for $|n|>1$ and $\t>\t_{0}+1$,
$$ | n| \ge \NN_{0} \= \left( { C_{0} \over  C\e_{0} }
\right)^{1/(\t_{0}+1)} ,
\tag 6.5 $$
%
%
%
%
%
We can define a map $t\to\e(t)$ such that
$$f_{n,m}(\e(t))=\O n-\sqrt{\o_m^2+n\n_{n,m}(\e)}   = {2\g t\over |n|^\t},\quad t\in [-1,1]\tag 6.6 $$ 
describes the interval defined by 6.2; then one has
$$\int_{\II_1} d\e=
\sum_{|n|\ge \NN_0, 0\le m\le 
c_2 \sqrt{n}
} \int_{-1}^{1} {\text d}t \left|
{{\text d}\e(t) \over {\text d}t} \right|,
\tag 6.7 $$
We have from the definition of $f_{n,m}$:
$$ {{\text d} f_{n,m} \over {\text d}t} = {{\text d} f_{n,m} \over {\text d} \e}
{{\text d} \e \over {\text d} t} = { 2\g \over |n|^{\t} } ,
\tag 6.8 $$
We need a lower bound on 
$$|{{\text d} f_{n,m}(\e)\over {\text d} \e}|= 
|n +n {\partial_\e \n_{n,m}\over 2\sqrt
{\o^2_m+n\n_{m}}}|\ge |n|- {C |n|\over \o_m}\tag 6.9$$
By 6.5 and the fact that, by 6.3,
${|n|\over \o_m}\le \bar C$
we get for $\e_0$ small enough
$$|{{\text d} f_{n,m}(\e)\over {\text d} \e}|>{|n|\over 2}\tag 6.10$$
We substitute in 6.7: 
$$\int_{\II_1} {\text d}\e \le\sum_{|n|\ge\NN_0,0<m\le c_2|n|^{1/2}}{C_1\over |n|^{\t+1}}
\le C_2\e_0^{(\t-{1\over 2}){1\over\t_0+1}}.$$
So the Cantor set of the $\e$ verifying 1.25 has relative measure 
$\to 1 $ as $\e_0\to 0$ if $\t>\t_0+{5\over 2}$.
\smallskip 
{\bf 6.3)} We define 
$\II_2$ the set of $\e\in (0,\e_0)$ belonging to $\II_1$ and verifying, 
for $m_1\not=m_2$, the condition $|\o_1 |n_i|-m_i^2|\le 1+\e_0 |n_i|$,
$i=1,2$, if $n=n_2-n_1$

$$|\O (n_2-n_1)\pm \sqrt{\o_{m_1}^2+n_1\n_{n_1,m_1}(\e)}\mp  
\sqrt{\o_{m_2}^2+(n_2)\n_{n_2,m_2}(\e)}|\le
{2\g\over |n_2-n_1|^{\t_0}}\tag 6.11$$
Of course $|n_i|\le C_1 m_i^2$, for $i=1,2$. 

For simplicity we choose the signs in 6.11 as in 
$-,+$ in (the other case
is done in the same way); then 6.11 can be verified for some $\e$
only if $m_1>m_2$. It holds that 
$$m_1^2-m_2^2\le
(\o_1+\e_0)|n|+1\tag 6.12$$
The proof is by contradiction; if it is not true
then  $m_1^2-m_2^2>
(\o_1+\e_0)|n_1-n_2|+1$ which implies $|n_1-n_2|\not=[\O^{-1}|m_1^2-m_2^2|]$,
where $[...]$ denotes the closest integer. Then
$$|\O n-\sqrt{\o_{m_1}^2+n_1\n_{n_1,m_1}}+
\sqrt{\o_{m_2}^2+n_2\n_{n_2,m_2}}|\ge  $$ 
$$|\O n-m_1^2+m_2^2 |
-{\m_0 \over m_1^2}-{\m_0\over m_2^2}- ({n_1\over m_1^2}+{n_2\over m_2^2})c_1\e_0|\ge |\O n- m_1^2+m_2^2 |- {1\over 4}(1+C\e_0)
\ge {1\over 8} \tag 6.13$$ 
as $|n_i|\le C_1 m_i^2$, for $i=1,2$, in contradiction with 6.11.

Then by 6.12 we get
$m_1+m_2\le C_2 {|n|\over m_1-m_2}\le C_2 |n|$
as $m_1-m_2\ge 1$;
hence $m_1\le C_3 |n|$ and $m_2\le C_3 |n|$.

Finally when the conditions 6.11 are satisfied, one has, for 
$n_1+n_2=n$ and $C_0=2\g$
$$ \eqalign{
2C_{0} | n|^{-\t_{0}} & \le
| \o_{1} n - (\o_{m_2}-\o_{m_1}) | \cr
& \le | \o_{1} n + \e n -
\sqrt{(\o_{m_2} + n_2 \n_{n_2,m_2}}+
\sqrt{(\o_{m_1} + n_1 \n_{n_1,m_1}}|
\cr
& +|\sqrt{(\o_{m_1} + n_1 \n_{n_1,m_1}}-\o_{m_1}|
+|\sqrt{(\o_{m_2} + n_2 \n_{n_2,m_2}}-\o_{m_2}|
+\e_0 |n|\cr
& \le C_{0} | n|^{-\t} + \e_{0}| n|+{|n_1||\n_{n_1,m_1}|\over m_1^2}+
{|n_2||\n_{n_2,m_2}|\over m_2^2} , \cr}
\tag 6.14 $$
now as $|n_i|\le c_2 m_i^2$ we have  that
$$ | n| \ge \NN_{1} \= \left(C_8 \e_{0} 
\right)^{1/(\t_{0}+1)} .
\tag 6.15 $$
We define the map $t\to \e(t)$ implicitly by:
$$f_{n,n_1,m_1,m_2}\equiv
\O n  -\sqrt{\o_{m_1}^2 +n_1\n_{m_1}}+ 
\sqrt{\o_{m_2}^2 +n_2\n_{m_2}}= {2\g t \over |n|^\t} \tag 6.16$$ 
We write
$$ \int_{\II_2} {\text d}\e =
\sum_{
n\le \left(C_8 \e_{0} 
\right)^{1/(\t_{0}+1)}; m_1, m_2\le C_3 |n|;|n_1|\le C m_1^2} 
\int_{-1}^{1} {\text d}t \left|
{{\text d}\e(t) \over {\text d}t} \right|\tag 6.17$$
We need a lower bound on 
$$|{{\text d} f_{n,n_1,m_1,m_2}(\e)\over {\text d} \e}|= 
|n-{n_1\partial_\e\n_{n_1,m_1}\over 
2\sqrt{\o_{m_1}^2+n_1\n_{m_1}}}+{n_2 \partial_\e\n_{m_2}
\over 2\sqrt{\o_{m_2}^2+n_2 \n_{nm_2}}}| \ge$$ 
$$ |n|- {(C+\e_0) |n_1|\over \o_{m_1}}- {(C+\e_0) 
|n_2|\over \o_{m_2}}\ge {|n|\over 2}\tag 6.18$$
where we have used that
${|n_i|\over \o_{m_i}}$ is bounded by a constant 
,(22), and we have chosen $\e_0$ small enough.
Hence we get
$$ \int_{\II_2} {\text d}\e =
\sum_{
(|n|\le \left(C_8 \e_{0} 
\right)^{1/(\t_{0}+1)}} C |n|^{-\t-1+4} 
\le \e_0^{-\t+4\over\t_0+1}\tag 6.19$$
so the Cantor set of the $\e$ verifying 1.26 has relative measure 
$\to 1 $ as $\e_0\to 0$ if $\t>\t_0+{5\over 2}$.
Finally we define $\II_3$ the set of $\e\in (0,\e_0)$ verifying 
$$|\O n\pm \sqrt{\o_{m_1}^2+n_1\n_{m_1}(\e)}\pm  
\sqrt{\o_{m_2}^2+(n+n_1)\n_{m_2}(\e)}|\le
{2\g\over |n|^{\t_0}}\tag 6.20$$
and one proceeds as above with the only difference that
6.20 can be true only if $|m_1^2+m_2^2|\le C_2 |n|$
hence $m_i\le C_2 \sqrt{|n|}$, $i=1,2$.
\smallskip\smallskip
\bigskip
\noindent{\bf Appendix A1. Measure of the set $M(\g)$}
\smallskip\smallskip
The analysis is very similar to the one in \S 6.
We call $J_1$ the set of $\m$ which do not satisfy the first condition in (kkl).
 $J_1$ is given by:
$$(1+\m)|n|-\sqrt{m^4+\m}= t{\g\over |n|^{\t_0}}\,,\quad t\in (-1,1);\tag A1.1 $$
the left hand side can be smaller than $1$ 
only if $n=[{\sqrt{m^4+\m}\over (1+\m)}] $, where $[...]$ is the closest integer; 
this implies that $m<c \sqrt{|n|}$ for a suitable constant $c$. Then
A1.1 defines the values $\m=\m(t)$ in $J_1$ so that:
$${\text meas}(J_1)= \sum_{n\atop m<c \sqrt{|n|}}\int_{-1}^1 dt|{d\m(t)\over dt}| 
\le \sum_{n\atop m<c\sqrt{|n|}}{2\g\over |n|^{\t_0+1}}\le 2\g\tag A1.2$$ 
as $|\partial_\m [(1+\m)|n|-\sqrt{m^4+\m}]|\ge {|n|\over 2}$, and $\t_0>{1\over 2}.$

 Let us call $J_2$ the set of $\m$ such that
$$(1+\m)|n|-\sqrt{m_1^4+\m}+ \sqrt{m_2^4+\m}=t{\g\over |n|^{\t_0}}\,,\quad t\in (-1,1);\tag A1.3$$
the left hand side can be smaller than one only if 
$n=[{\sqrt{m_1^4+\m}-\sqrt{m_2^4+\m}\over (1+\m)}] $ which implies $|m_1^2-m_2^2|<c|n|$ and therefore $m_1+m_2<c|n|$.
$${\text meas}(J_2)= \sum_{n\atop m_1,m_2<c|n|}\int_{-1}^1 dt|{d\m(t)\over dt}| 
\le \sum_{n}{2\g\over |n|^{\t_0-1}}\le 2\g, \tag A1.4$$ 
as $\t_0>2$.
Finally we proceed in the same way for $J_3$ the set of $\m$ which do not satisfy the third condition in (kkl). We have proved that the complementary set to $M(\g)$ is of order $6\g<{1\over 8}$ provided that $\g$ is small enough,
that is $\g\le 2^{-6}$.
\smallskip
\bigskip
\noindent{\bf Appendix A2. Proof of 1.7}
\bigskip
The equation for the coefficients 1.10 follows immediately from
$$\int_0^\pi d x \sin(m x) \sin(m_1 x)=\pi \d_{m,m_1}\tag A2.1$$
and 
$$\int_0^\pi d x \sin(m x)\sin(m_1 x)\sin(m_2 x)
=\sum_{\e,\e_1,\e_2=\pm} (\e\e_1\e_2) 
{e^{i(\e m+\e_1 m_1+\e_2 m_2)\pi}-1\over i
(\e m+\e_1 m_1+\e_2 m_2)}\tag A2.2$$
which is vanishing if $\pm m_1\pm m_1\pm m_2$ is even, while
it it is odd it is equal to
$$4[{1\over m+m_1+m_2}-{1\over m+m_1-m_2}-{1\over m-m_1+m_2}+
{1\over m-m_1-m_2}]=$$
$${8 m_1 m_2 m\over (m^2-(m_1-m_2)^2) 
(m^2-(m_1+m_2)^2)
}\tag A2.3$$

\smallskip
\bigskip
\noindent
{\bf Appendix A3. Proof of the Lemmas 5 and 8}
\smallskip 
In order to prove Lemma 5 
we prove inductively the bound, for $\th\in\Th_{R,n.m}$
$$ N^*_{h}(\th) \le \max\{ 0 , 2 K(\th) 2^{(2-h)/\t} - 1 \} ,
\tag A3.1 $$
where $N^*_{h}(\th)$ is the number of  non resonant lines.
As we are supposing $\Val(\th)\not=0$ it holds
for any $\ell$ that
$\g 2^{-h_\ell-1}\leq |x_{n_\ell,m_\ell}|\leq \g 2^{-h_\ell+1}$.
This implies, by the first Diophantine condition, that $\th$ can have a line on scale
$h$ only  if $K(\th) > 2^{(h-1)/\t}$.
Then  one can have $N_{h}(\th) \ge 1$ only
if $K(\th)$ is such that $K(\th) > k_{0}\=2^{(h-1)/\t}$:
therefore for values $K(\th) \le k_{0}$ the bound (4.9) is satisfied.

If $K(\th)>k_{0}$, we assume that
the bound holds for all trees $\th'$ with $K(\th')<K(\th)$.
Define $E_{h}=2^{-1}(2^{(4-h)/\t})^{-1}$: so we have
to prove that $N^*_{h}(\th)\le \max\{0,K(\th) E_{h}^{-1}-1\}$.

Call $\ell$ the root line of $\th$ and
$\ell_{1},\ldots,\ell_{m}$ the $m\ge 0$ lines on scale $\ge h$
which are the closest to $\ell$ 

If the root line $\ell$ of $\th$ is  on scale $< h$
 then
$$ N_{h}^*(\th) = \sum_{i=1}^{m} N_{h}^*(\th_{i}) ,
\tag A3.2 $$
where $\th_{i}$ is the subtree with $\ell_{i}$ as root line,
hence the bound follows by the inductive hypothesis.

If the root line $\ell$ has scale $\ge h$
then $\ell_{1},\ldots,\ell_{m}$ are the entering line of a cluster $T$.

By denoting again with $\th_{i}$ the subtree
having $\ell_{i}$ as root line, one has
$$ N_{h}^*(\th) = 1 + \sum_{i=1}^{m} N^*_{h}(\th_{i}) ,
\tag A3.3$$
and the bound becomes trivial if either $m=0$ or $m\ge 2$.

If $m=1$ then one has a cluster $T$ with two external lines
$\ell$ and $\ell_{1}$, with $h_{\ell_1}, h_{\ell}\ge h $ so that by 1.14:
$$ \left| |\O n_{\ell}| - \sqrt{\o_{m_{\ell}}+n\n_{m_\ell,n_\ell}} \right|
\le 2^{-h+1} \g , \qquad
\left| |\O n_{\ell_{1}}| - \sqrt{\o_{m_{\ell_{1}}}+n\n_{m_{\ell_1},n_{\ell_1}}} \right|
\le 2^{-h+1} \g ,
\tag A3.4 $$
As $\ell$ is non resonant, surely  $n_{\ell} \not = n_{\ell_{1}}$ (otherwise
if $n_{\ell} = n_{\ell_{1}}$ then $m_{\ell} \not = m_{\ell_{1}}$
hence the two lines cannot have both scale $\ge h$).
Hence by 1.26
one has, for suitable $\h_{\ell},\h_{\ell_{1}}\in\{+,-\}$,
$$ 2^{-h+2} \g \ge \big| \O
( n_{\ell}-n_{\ell_{1}}) + \h_{\ell}\sqrt{\o_{m_{\ell}}+n\n_{m_\ell,n_\ell}} +
\h_{\ell_{1}}\sqrt{\o_{m_{\ell_{1}}}+n\n_{m_{\ell_1},n_{\ell_1}}}  \big| \ge
\g |n_{\ell}-n_{\ell_{1}}|^{-\t} ,
\tag A3.5$$
so that $K(\th)-K(\th_{1})> E_{h}$.
 Hence 
by using the inductive hypothesis
$$ \eqalign{
N^*_{h}(\th) & = 1 + N^*_{h}(\th_{1})
\le 1 + K(\th_{1}) E_{h}^{-1} - 1 \cr
& \le 1 + \Big( K(\th) - E_{h}\Big) E_{h}^{-1} - 1
\le K(\th) E_{h}^{-1} - 1 , \cr}
\tag A3.6$$
hence the bound is proved also
if the root line is on scale $\ge h$.
\smallskip 
We prove Lemma 8 for $\LL\Val(\hat\th)$, 
$\hat\th\in\R_{R,n,m,h}$. We consider the two subtrees
entering in $v_0$; one, called $\tilde\th$, does not contain the endpoint $e$
and the bounds of the preceding lemma can be applied, so
we consider the subtree $\th$ containing $e$.
We proceed inductively on $h$ for $\th$ proving that
$N_{h}^*(\th)\leq 2K(\th)2^{2-h\over \t}$; such bound and the bound A3.1 for $\tilde\th$
immediately implies 4.7 as $K(\hat\th)=K(\tilde\th)+K(\th)-1$.

In order to prove $N_{h}^*(\th)\leq 2K(\th)2^{2-h\over \t}$
we define $k_{0}=2^{(h-1)/\t}$. One has $N_{h}^{*}(\th)=0$ for
$K(\th)<k_{0}$. In fact if 
the line $\ell$ with scale $h$ is not in the path, 
this follows from the first Melnikov condition.
If such line is 
on the path we have that, if $\LL\Val(\th)$ is non vanishing
$$||\O n^0_\ell+\bar\o_{n,m}|- 
\sqrt{\o_{m_\ell}+n\n_{m_\ell,n_\ell }})|\leq \g 2^{-h+1}\tag A3.7$$
so that by the second Melnikov condition:
$$K(\th)\ge |n^0_\ell|\ge 2^{h-1\over \t}.\tag A3.8$$

Then for $K(\th)< k_0$ the bound is satisfied;
for $K \ge k_{0}$, we assume that the bound 
holds for all $K(\th)=K'<K$,
and we show that it follows also for $K(\th)=K$.
If $K(\th)>k_{0}$, we assume that
the bound holds for all trees $\th'$ with $K(\th')<K(\th)$.
Define $E_{h}=2^{-1}(2^{(4-h)/\t})^{-1}$: so we have
to prove that $N^*_{h}(\th)\le K(\th) E_{h}^{-1}$.
If the root line $\ell$ of $\th$ is  on scale $< h$
then
$$ N_{h}^*(\th) = \sum_{i=1}^{m} N_{h}^*(\th_{i}) ,
\tag A3.9$$
where $\th_{i}$ is the subtree with $\ell_{i}$ as root line,
hence the bound follows by the inductive hypothesis.
If the root line $\ell$ has scale $\ge h$
then $\ell_{1},\ldots,\ell_{m}$ are the entering line of a cluster $T$.
The same occurs if the root line is on scale $\ge h$ and non-resonant,
and, by calling $\ell_{1},\ldots,\ell_{m}$ the lines on scale
$\ge h$ which are the closest to $\ell$, one has $m \ge 2$:
in fact in such a case at least $m-1$ among the subtrees
$\th_{1},\ldots,\th_{m}$ verifies the bound A3.1 so that
$$ \eqalign{
N_{h}^{*}(\th) & = 1 + \sum_{i=1}^{m} N_{h}^{*}(\th_{i}) \cr
& \le 1 + E_{h}^{-1} \sum_{i=1}^{m} K(\th_{i}) - (m-1)
\le E_{h}  K(\th) , \cr}
\tag A3.10$$
If $m=0$ then $N_{h}^{*}(\th)=1$ and $K(\th) 2^{(2-h)/\t} \ge 1$
because one must have $K(\th) \ge k_{0}$.
So the only non-trivial case is when
one has $m=1$.
In this case $\ell,\ell_{1}$ are on the path 
connecting the external lines of the resonance
$$||\O n^0_\ell+\bar\o_{n,m}|- 
\sqrt{\o_{m_\ell}+n\n_{m_\ell,n_\ell }})|\leq \g 2^{-h+1}\tag A3.11$$
$$||\O n^0_{\ell_1}+\bar\o_{n,m}|- 
\sqrt{\o_{m_{\ell_1}}+n\n_{m_{\ell_1},n_{\ell_1}}})|\leq \g 2^{-h+1}\tag A3.12$$
so that for suitable $\h_{\ell},\h_{\ell_{1}}\in\{+,-\}$
$$ 2^{-h+2} \g \ge \big| \O
(n^0_{\ell}-n^0_{\ell_{1}}) + \h_{\ell}\sqrt{\o_{m_{\ell}}+
n\n_{m_{\ell},n_{\ell}}}+
\h_{\ell_{1}} \sqrt{\o_{m_{\ell_1}}+n\n_{m_{\ell_1},n_{\ell_1}}} \big| \ge
\g |n^0_{\ell}-n^0_{\ell_{1}}|^{-\t} ,
\tag A3.13$$
from which $K(\th)-K(\th_{1})\ge |n^0_{\ell}-n^0_{\ell_{1}}|\ge 
 2^{(h-2)/\t}$ and by the analogous of A3.6
the bound is proved.

\bigskip
\noindent
{\bf Appendix A4. Proof of Lemma 7}
\bigskip
In order to prove  4.15 we proceed by induction; 
consider a tree $\th$ with $k$ internal nodes and $s_v=2$ for any $v$; 
we call $\ell_1,\ell_2$ the two lines  entering   
$v_0$; we call $m_{\ell_1}=m_1$ and $m_{\ell_2}=m_2$
the root lines
of two subtrees $\th_1$ and $\th_2$ with $k_1\geq 0$ and $k_2\geq 0$ vertices,
and $k_1+k_2=k-1$.
If  $k_1=0$ (or $k_2=0$) then one of the two cases holds:

1. $\ell_1$ connects to an end-node so that $|\h_{v}|=1 $ and $m_1=1$.

2. $\ell_1$ connects to the external node so that $|\h_{v}|=|\h_e|=1/m_e^3 $ (this case is possible only of $\th\in 
\RR_{R,h,n,m}$).

So we can proceed with our inductive hypothesis and suppose that our bound holds for all trees with $0\leq k_1<k$ end-nodes.
Without loss of generality we can suppose  that $m_1\ge m_2$. We perform the bound  
$m+m_1+m_2>m$ ,$m+m_1-m_2>m$  so that:
$$\sum^*_{  m_1\ge m_2} 
{1\over |(m^2-(m_1+m_2)^2)| |(m^2-(m_1-m_2)^2)|}
{1\over m_1^2m_2^2}\le$$ 
$${1\over m^2}\sum^*_{ 
m_1\ge m_2} {1\over m_1^2m_2^2}{1\over |m-m_1-m_2||m-m_1+m_2|}\tag A4.1$$ 
 If $m_1\le {m\over 4}$ then the bound is trivial as:
$${1\over m^2}\sum^*_{{m\over 4}\ge  m_1\ge m_2} 
{1\over |m-m_1+m_2| |m-m_1-m_2|}{1\over
m_1^2m_2^2}\le {8\over 3 m^4}\sum^*_{m_1,m_2}
{1\over m_1^2m_2^2}\le {C_1\over m^4}.\tag A4.2$$
 
In the remaining terms we treat separately the cases   $m_1\leq m-1$, $m_1>m$ and $m_1=m$.

We notice that in the first case $|m-m_1+m_2|\geq m-m_1$, while in the second case $|m-m_1-m_2|\geq m_1-m$ . 
We then obtain the bound:
 $${1\over m^2}\sum^*_{ m_1>{m\over 4}\,,\, 
m_1\ge m_2} {1\over m_1^2m_2^2}{1\over |m-m_1-m_2||m-m_1+m_2|}\leq $$ 
$${1\over m^2}
\Big(\sum^*_{m-1\ge m_1> {m\over 4}}
{1\over m_1^2(m- m_1)}(\sum_{m_2=1}^\io {1\over m_2^2|m-m_1-m_2|})+\tag A4.3$$ 
$$\sum^*_{m_1> m} {1\over m_1^2(m_1-m)}(\sum_{m_2=1}^\io 
{1\over m_2^2|m-m_1+m_2|})+ {1\over m^2}\sum_{m_2=1}^\io {1\over m_2^2} \Big)$$   
Now we estimate the sums with integrals:
$$\sum_{n\ne A}{1\over |A-n|n^2}\le C_0[
\int_{x=1}^{A-1}{1\over (A-x)x^2}+
\int_{x=A+1}^{\io}{1\over (x-A)x^2}+ {C_2\over A^2}]\le  
{C_3\over A}\tag A4.4$$
$$\sum_{n\neq A}{1\over (A-n)^2n^2}\le C_0[
\int_{x=1}^{A-1}{1\over (A-x)^2x^2}+\int_{x=A+1}^{\io}
{1\over (x-A)^2x^2}+ {C_2\over A^2}]\le {C_3\over A^2} \tag A4.5$$
We use the first bound on the sum over $m_2$ then 
in the sum over $m_1$ we obtain a series as in the second bound 
immediately implying
$$\sum^*_{ m_1>{m\over 4}\,,\, m_1\ge m_2} 
{1\over m_1^2m_2^2} {1\over |(m^2-(m_1+m_2)^2)| 
|(m^2-(m_1-m_2)^2)|}{1\over m_1^2m_2^2}\le {\bar C\over m^4}\tag A4.6).$$
This implies the inductive hypothesis, as $\th$
has $k$ vertices  and $k_1+k_2=k-1$, by choosing $ C_1=\bar C$ in  4.15.

\medskip

E-mail address :  mastropi\@mat.uniroma2.it; procesi\@mat.uniroma3.it

\medskip

Received April 2005

\enddocument

%

\rife{Bo1}{2}{
J. Bourgain,
{\it Construction of quasi-periodic solutions for Hamiltonian
perturbations of linear equations and applications to nonlinear PDE},
Internat. Math. Res. Notices
{\bf 1994}, no. 11, 475ff., approx. 21 pp. (electronic). }
\rife{Bo2}{3}{
J. Bourgain,
{\it Construction of periodic solutions of nonlinear
wave equations in higher dimension},
Geom. Funct. Anal.
{\bf 5} (1995), 629--639. }
\rife{Bo4}{4}{
J. Bourgain,
{\it Periodic solutions of nonlinear wave equations},
Harmonic analysis and partial differential equations
(Chicago, IL, 1996), 69--97, Chicago Lectures in Math.,
Univ. Chicago Press, Chicago, IL, 1999. }
\rife{Bo2}{5}{
J. Bourgain,
{\it Construction of periodic solutions of nonlinear
wave equations in higher dimension},
Geom. Funct. Anal.
{\bf 5} (1995), 629--639. }
\rife{CK}{6}{
M.Countryman and R.Kannan,
{\it Nonlinear Damped vibration of beams with periodic forcng},
Int. J. Non-Linear Mech. 27,1, (1992), 75--83. }
\rife{CW1}{7}{
W. Craig and C.E. Wayne,
{\it Newton's method and periodic solutions of nonlinear
wave equations},
Comm. Pure Appl. Math.
{\bf 46} (1993), 1409--1498. }
\rife{E}{8}{
L.H. Eliasson,
{\it Absolutely convergent series expansions for quasi periodic motions},
Math. Phys. Electron. J.
{\bf 2} (1996), Paper 4, 33 pp. (electronic). }
\rife{Ga}{9}{
G. Gallavotti,
{\it Twistless KAM tori},
Comm. Math. Phys.
{\bf 164} (1994), no. 1, 145--156. } 
\rife{GM1}{10}{
G. Gentile and V. Mastropietro,
{\it 
Convergence of Lindstedt series for the nonlinear wave equation}.  
Commun. Pure Appl. Anal.  3  (2004),  no. 3, 509--514}
\rife{GM2}{11}{
G. Gentile and V. Mastropietro,
{\it Construction of periodic solutions of the nonlinear wave equation
with Dirichlet boundary conditions by the Linsdtedt series method},
J. Math. Pures Appl. (9),  83  (2004),  no. 8, 1019--1065. }
\rife{GMP}{12}{
G. Gentile, V. Mastropietro and M.Procesi,
{\it  Periodic solutions for completely resonant nonlinear
wave equations
with Dirichlet boundary conditions},
Comm. Math. Phys.(2005)}
\rife{GoR}{13}{
C. Godsil, G. Royle,
{\it Algebraic graph theory},
Graduate Texts in Mathematics 207.
Spring\-er, New York, 2001. }
\rife{HP}{14}{
F. Harary, E.M. Palmer,
{\it Graphical enumeration},
Academic Press, New York-London, 1973. }
\rife{HP1}{15}{
W.T. van Horseen, {\it An asymptotic theory for a Class of Initial-Boundary
Value for weakly nonlinear Wave equations} SIAM J. Applied Math.,
48, (1988), 1227-1243}
\rife{HLM}{16}{
Hsu, Lu, Machino, {\it Periodic solutions to the 1-dimensional compressible Euler equation}, preprint}
\rife{K1}{17}{
S.B. Kuksin,
{\it Nearly integrable infinite-dimensional Hamiltonian systems},
Lecture Notes in Mathematics 1556, Springer, Berlin, 1994. }
\rife{KP}{18}{
S.B. Kuksin, J.P\"oschel,
{Invariant Cantor manifolds of quasi-periodic oscillations
for a nonlinear Schr\"odinger equation},
Ann. of Math. (2)
{\bf 143} (1996), no. 1, 149--179. }
\rife{L}{19}{
A. M. Lyapunov,
{\it Probl\`eme g\'en\'eral de la  stabilit\'e du mouvement},
Ann. Sc. Fac. Toulouse
{\bf 2} (1907), 203--474. }
\rife{Wh}{20}{
H. Whitney,
{\it Analytic extensions of differential
functions defined in closed sets},
Trans. Amer. Math. Soc.
{\bf 36} (1934), no. 1, 63--89. }

\biblio
\ciao